\documentclass[a4paper,11pt,english]{smfart}

\usepackage{amssymb}
\usepackage{url}
\usepackage[T1]{fontenc}
\RequirePackage{calrsfs}
\DeclareSymbolFont{rsfscript}{OMS}{rsfs}{m}{b}
\DeclareSymbolFontAlphabet{\mathrsfs}{rsfscript}
\usepackage[utopia,expert]{mathdesign}
\usepackage{lscape}
\usepackage{array, boldline, makecell, booktabs}

\usepackage{enumitem}

\usepackage{todonotes}

\usepackage[leftbars]{changebar}

\usepackage{palatino}
\usepackage{rotating}
\usepackage{graphicx}

\usepackage{tikz-cd}

\usepackage{enumerate}

\usepackage{amscd}

\usepackage{color}
\definecolor{shadecolor}{gray}{0.90}

\def\bfit{\bfseries\itshape}

\def\proofname{D\'emonstration}

\input xypic
\xyoption{all}
\xyoption{arc}

\makeindex

\def\indexname{Index des notations}

\newtheorem{theo}{Theorem}[section]
\def\thetheo{\thesection.\arabic{theo}}
\newtheorem{prop}[theo]{Proposition}

\newtheorem{lem}[theo]{Lemma}

\renewcommand\theequation{\thetheo}
\def\equat{\refstepcounter{theo}\begin{equation}}
\def\endequat{\end{equation}}

\renewcommand\thetable{\Roman{table}}
\renewcommand\thefigure{\Roman{figure}}

\renewcommand\thesection{\arabic{section}}
\renewcommand\thesubsection{\thesection.\Alph{subsection}}

\def\contentsname{Table des mati\`eres}
\def\listfigurename{Liste des figures}
\def\listtablename{Liste des tables}
\def\appendixname{Appendice}

\def\refname{R\'ef\'erences}

    \def\CM{{\mathbb{C}}}

    \def\FM{{\mathbb{F}}}

    \def\QM{{\mathbb{Q}}}
    \def\RM{{\mathbb{R}}}
\def\SG{{\mathfrak S}}

\def\Ab{{\mathbf A}}    
    
    \def\CC{{\mathcal{C}}}
    
    \def\EC{{\mathcal{E}}}
    
\def\Gb{{\mathbf G}}

\def\Lb{{\mathbf L}}    \def\LC{{\mathcal{L}}}

\def\Ob{{\mathbf O}}    
\def\Pb{{\mathbf P}}

\def\Sb{{\mathbf S}}    \def\SC{{\mathcal{S}}}
    
    \def\UC{{\mathcal{U}}}

    \def\XC{{\mathcal{X}}}
    
    \def\ZC{{\mathcal{Z}}}

\def\Nrm{{\mathrm{N}}}

\def\Srm{{\mathrm{S}}}

\def\Wrm{{\mathrm{W}}}

\def\Zrm{{\mathrm{Z}}}

\def\g{\gamma}
\def\G{\Gamma}
\def\d{\delta}
\def\D{\Delta}
\def\e{\varepsilon}
\def\ph{\varphi}
\def\l{\lambda}
\def\L{\Lambda}

\def\O{\Omega}
\def\r{\rho}
\def\s{\sigma}

\def\t{\tau}

\def\z{\zeta}

\def\mub{{\boldsymbol{\mu}}}

\DeclareMathOperator{\Gal}{{\mathrm{Gal}}}

\def\to{\rightarrow}

\def\injto{\hookrightarrow}

\def\DS{\displaystyle}

\def\finl{~$\blacksquare$}

\def\lexp#1#2{\kern\scriptspace\vphantom{#2}^{#1}\kern-\scriptspace#2}
\def\le{\hspace{0.1em}\mathop{\leqslant}\nolimits\hspace{0.1em}}
\def\ge{\hspace{0.1em}\mathop{\geqslant}\nolimits\hspace{0.1em}}

\mathchardef\inferieur="321E
\mathchardef\superieur="321F

\def\eqna{\begin{eqnarray*}}
\def\endeqna{\end{eqnarray*}}

\def\itemth#1{\item[${\mathrm{(#1)}}$]}

\catcode`\@=11
\long\def\@car#1#2\@nil{#1}
\long\def\@first#1#2{#1}
\long\def\@second#1#2{#2}
\long\def\ifempty#1{\expandafter\ifx\@car#1@\@nil @\@empty
  \expandafter\@first\else\expandafter\@second\fi}
\catcode`\@=12

\def\GL{{\mathrm{GL}}}

\DeclareMathOperator{\Ref}{Ref}

\def\boitegrise#1#2{\begin{centerline}{\fcolorbox{black}{shadecolor}{~
    \begin{minipage}[t]{#2}{\vphantom{~}#1\vphantom{$A_{\DS{A_A}}$}}
            \end{minipage}~}}\end{centerline}\medskip}

\theoremstyle{remark}
\newtheorem{rema}[theo]{Remark}

\newtheorem{exemple}[theo]{Example}
\newtheorem{exemples}[theo]{Examples}

\theoremstyle{plain}

\def\BIL{LR}
\def\GAUCHE{L}
\def\CAR{CAR}
\def\FAM{FAM}

\def\xyinj{\ar@{^{(}->}}
\def\xysur{\ar@{->>}}

\bigskip

\def\petitespace{\vphantom{$\DS{\frac{\DS{A^A}}{\DS{A_A}}}$}}

\makeatletter
\def\hlinewd#1{%
\noalign{\ifnum0=`}\fi\hrule \@height #1 %
\futurelet\reserved@a\@xhline}
\makeatother

\newlength\epaisLigne

\usepackage{multirow,multicol}

\makeatletter
\def\hlinewd#1{%
\noalign{\ifnum0=`}\fi\hrule \@height #1 %
\futurelet\reserved@a\@xhline}
\makeatother

\usepackage{multirow}

\def\petitespace{\vphantom{$\DS{\frac{\DS{A^A}}{\DS{A_A}}}$}}

\usepackage{lscape}

\def\GL{\operatorname{\Gb\Lb}\nolimits}

\def\aff{{\mathrm{aff}}}
\def\sing{{\mathrm{sing}}}
\def\irr{{\mathrm{irr}}}

\begin{document}

\title{Some singular curves and surfaces arising \\ 
from invariants of complex reflection groups}

\author{{\sc C\'edric Bonnaf\'e}}
\address{C. Bonnaf\'e, IMAG, Universit\'e de Montpellier, CNRS, Montpellier, France} 

\makeatletter
\email{cedric.bonnafe@umontpellier.fr}
%
%

\date{\today}

\thanks{The author is partly supported by the ANR 
(Project No ANR-16-CE40-0010-01 GeRepMod)}

\begin{abstract}
We construct highly singular projective curves and surfaces defined 
by invariants of primitive complex reflection groups.
\end{abstract}

\maketitle
\pagestyle{myheadings}
\markboth{\sc C. Bonnaf\'e}{Singular curves and surfaces}

It is a classical problem to determine the maximal number of singularities 
of a given type that a curve or a surface might have. Several kinds of upper 
bounds have been given~\cite{sakai},~\cite{bruce},~\cite{miyaoka},~\cite{varchenko},~\cite{wahl}..., 
and these bounds have been approached 
for small 
degrees~\cite{ivinskis},~\cite{barth},~\cite{escudero},~\cite{escudero 2},~\cite{endrass},~\cite{triple},~\cite{labs},
~\cite{sarti},~\cite{stagnaro}... 
or general degrees~\cite{chmutov}. 

In~\cite{barth},~\cite{sarti 0},~\cite{sarti 1},~\cite{sarti}, 
Barth and Sarti used pencils of surfaces constructed from invariants 
of some finite Coxeter subgroups of $\Gb\Lb_4(\RM)$ to obtain surfaces of degree $6$, $10$, $12$ 
with the biggest number of nodes known up to now. We have decided 
to explore more systematically pencils of curves and surfaces constructed from invariants 
of finite complex reflection subgroups of $\Gb\Lb_3(\CM)$ or $\Gb\Lb_4(\CM)$. In this paper, we gather 
the results of these computations (made with {\sc Magma}~\cite{magma}) obtained from the 
{\it primitive} complex reflection. As the reader will see, not all the primitive complex 
reflection groups lead to interesting examples but these investigations have lead 
to the discovery of the following curves or surfaces, which improve some known lower bounds and are 
quite close to upper bounds found by Sakai~\cite{sakai} for curves or Miyaoka~\cite{miyaoka} for surfaces 
(we refer to Shephard-Todd notation~\cite{shephard todd} for complex reflection groups; 
for Coxeter groups, we also use the notation $\Wrm(\G)$, where $\G$ is a Coxeter graph):
\begin{itemize}
\itemth{a} Using the complex reflection group $G_{24} \subset \Gb\Lb_3(\CM)$, we construct 
a curve of degree $14$ with $42$ cusps (i.e. singularities of type $A_2$): this improves 
known lower bounds (see Example~\ref{ex:record?}). 
Note that the known upper bound for the number of cusps 
of a curve of degree $14$ in $\Pb^2(\CM)$ is $55$.

\itemth{b} Using the complex reflection group $G_{26}$, we construct a curve of degree 
$18$ with $36$ singular points of type $E_6$ (see Example~\ref{ex:g26}). 
We do not know if such a bound was 
already reached. 

\itemth{c} Let $\mu_{D_4}(d)$ denote the maximal number of quotient singularities of 
type $D_4$ that an irreducible surface in $\Pb^3(\CM)$ might have. Miyaoka~\cite{miyaoka} proved that
$$\mu_{D_4}(d) \le \frac{16}{117}d(d-1)^2.$$
For $d=8$, $12$, or $24$, this reads
$$\mu_{D_4}(8) \le 53,\qquad \mu_{D_4}(12) \le 198\qquad\text{and}\qquad \mu_{D_4}(24) \le 1~\!736.$$
Using respectively the complex reflection groups $G_{28}=\Wrm(F_4)$, $G_{29}$ and $G_{32}$, 
we prove that 
$$\mu_{D_4}(8) \ge 48,\qquad \mu_{D_4}(12) \ge 160
\qquad\text{and}\qquad \mu_{D_4}(24) \ge 1~\!440
$$
(see Examples~\ref{ex:f4} and~\ref{ex:g29}(3) and Table~\ref{tab:singular surfaces}).
This improves considerably the last known lower bounds~\cite{escudero 2}. 
Recall that, by standard arguments, this implies that 
$\mu_{D_4}(8k) \ge 48 k^3$, $\mu_{D_4}(12k) \ge 160 k^3$ and $\mu_{D_4}(24k) \ge 1~\!440 k^3$ 
for all $k \ge 1$. Note that the fact that $\mu_{D_4}(24) \ge 1~\!440$ 
was first announced in~\cite{1440} (and a previous lower bound $\mu_{D_4}(8) \ge 44$ was 
also obtained): see Section~\ref{sec:g32} for details. 
\end{itemize}

\bigskip

We also found examples which do not improve known lower bounds but might possibly 
be interesting for the number and the type of singularities they contain (with 
``big'' multiplicities or ``big'' Milnor numbers): see Examples~\ref{ex:g26},~\ref{ex:g29},~\ref{ex:g31}. 
The examples might also be interesting for their big group of automorphisms.

These computations also show that Miyaoka bounds 
are quite sharp, even for singularities that are not of type $A$. 
Contrary to previous constructions, the singular points of 
our curves or surfaces are in general not all real\footnote{There is an important 
exception to this remark: all the singular points of the surface of degree $8$ with 
$48$ singularities of type $D_4$ constructed in Example~\ref{ex:f4} 
have rational coordinates.} 
(even though most 
of these varieties are defined over $\QM$). By contrast, note also that, 
using a theorem of Marin-Michel on automorphisms of reflection 
groups~\cite{marin-michel}, we can show that the 
Sarti dodecic can be defined over $\QM$ (this was still 
an open question). 

For the smoothness of the exposition, we have decided to include most of the {\sc Magma} codes 
in separate texts~\cite{1440} (for varieties associated with $G_{32}$) 
and~\cite{bonnafe calculs} (for the other examples), as well as some explicit polynomials: 
these two texts are not intended to be published, but are made for the reader 
interested in checking the computations by himself.

\bigskip

\section{Notation, preliminaries}

\medskip

We fix an $n$-dimensional $\CM$-vector space $V$ 
and a finite subgroup $W$ of $\Gb\Lb_\CM(V)$. We set 
$$\Ref(W)=\{s \in W~|~\dim_\CM(V^s)=n-1\}.$$

\bigskip

\boitegrise{{\bf Hypothesis.} {\it We assume throughout this paper that 
$$W=\langle \Ref(W) \rangle.$$
In other words, $W$ is a {\bfit complex reflection group}. We also assume that 
$W$ acts {\bfit irreducibly} on $V$. The number $n$ is called the {\bfit rank} of $W$.}}{0.75\textwidth}

\bigskip
\def\propor{\asymp}

\subsection{Invariants} 
We denote by $\CM[V]$ the ring of polynomial functions on $V$ (identified with the symmetric 
algebra $\Srm(V^*)$ of the dual $V^*$ of $V$) and by $\CM[V]^W$ 
the ring of $W$-invariant elements of $\CM[V]$. 
By Shephard-Todd/Chevalley-Serre Theorem~\cite[Theorem~4.1]{broue}, 
there exist $n$ algebraically independent 
homogeneous elements $f_1$, $f_2$,\dots, $f_n$ of $\CM[V]^W$ such that 
$$\CM[V]^W=\CM[f_1,f_2,\dots,f_n].$$
Let $d_i=\deg(f_i)$. We will assume that $d_1 \le d_2 \le \cdots \le d_n$. 
A family $(f_1,f_2,\dots,f_n)$ satisfying the above property 
is called a {\it family of fundamental invariants} of $W$. 
Whereas such a family is not uniquely defined, the list 
$(d_1,d_2,\dots,d_n)$ is well-defined and is called the list of {\it degrees} 
of $W$. If $f \in \CM[V]$ is homogeneous, we will denote by 
$\ZC(f)$ the projective (possibly reduced) hypersurface in $\Pb(V) \simeq \Pb^{n-1}(\CM)$  
defined by $f$. Its singular locus will be denoted by $\ZC_\sing(f)$. 
A homogeneous element $f \in \CM[V]Ŵ$ is called a {\it fundamental invariant} 
if it belongs to a family of fundamental invariants.

Recall that a subgroup $G$ of $\Gb\Lb_\CM(V)$ is called {\it primitive} if 
there does not exist a decomposition $V=V_1 \oplus \cdots \oplus V_r$ with $V_i \neq 0$ and 
$r \ge 2$ such that $G$ permutes the $V_i$'s. 
We will be mainly interested in {\it primitive} (often called {\it exceptional}) complex reflection groups, 
and we will refer to Shephard-Todd numbering~\cite{shephard todd} for such groups (there are $34$ 
isomorphism classes, named $G_i$ for $4 \le i \le 37$). Almost all the 
computations\footnote{Some Milnor and Tjurina numbers were computed with {\sc Singular}~\cite{singular}.}
have been done using the software {\sc Magma}~\cite{magma}. 

\bigskip
\def\QMov{{\overline{\QM}}}

\subsection{Marin-Michel Theorem} 
Let $\QMov$ denote the algebraic closure of $\QM$ in $\CM$ and we set 
$\G=\Gal(\QMov/\QM)$. 
Using the classification of finite reflection groups, 
Marin-Michel~\cite{marin-michel} proved that there exists a $\QM$-structure 
$V_\QM$ of $V$ such that:
\begin{itemize}
\itemth{1} $V_\QMov = \QMov \otimes_\QM V_\QM$ is stable under the action of $W$ 
(so that $W$ might be viewed as a subgroup of $\GL_\QMov(V_\QMov)$).

\itemth{2} The action of $\G$ on $\GL_\QMov(V_\QMov)$ induced by the 
$\QM$-form $V_\QM$ stabilizes $W$.
\end{itemize}
This implies that $\QMov[V_\QMov]$ is a $\QMov$-form of $\CM[V]$ 
stable under the action of $W$ and that the action of $\G$ on $\QMov[V_\QMov]$ 
induced by the $\QM$-form $\QM[V_\QM]$ stabilizes the invariant ring 
$\QMov[V_\QMov]$.

\medskip

\begin{prop}\label{prop:q}
The Sarti dodecic can be defined over $\QM$.
\end{prop}

\smallskip

\begin{rema}
An explicit polynomial with rational coefficients 
defining the Sarti dodecic is given in~\cite{bonnafe calculs}.\finl
\end{rema}

\smallskip

\begin{proof}
Assume here that $W$ is a Coxeter group of type $H_4$ acting on 
a vector space $V$ of dimension $4$. We fix a $\QM$-form $V_\QM$ 
as above. Let $f$ be a homogeneous invariant of $W$ of degree $12$ defining 
the Sarti dodecic: it belongs to $\QMov[V_\QMov]$. We fix a $\QM$-basis 
$(h_1,h_2,\dots,h_{455})$ of the homogeneous component of degree $12$ of 
$\QM[V_\QM]$. It is also a $\QMov$-basis of the homogeneous component of degree $12$ of 
$\QMov[V_\QMov]$. By multiplying $f$ by a scalar if necessary, 
we may assume that there exists $i \in \{1,2,\dots,455\}$ 
such that the coefficient of $f$ on $h_i$ is $1$. 

Now, if $\g \in \G$, then 
$\lexp{\g}{f}$ is also an invariant of $W$ of degree $12$ defining 
an irreducible projective surface with $600$ nodes. By the unicity 
of such an invariant~\cite{sarti}, this forces $\lexp{\g}{f}= \xi f$ 
for some $\xi \in \QMov^\times$. But $\xi=1$ because 
the coefficient of $f$ on $h_i$ is $1$. So $f \in \QM[V_\QM]$.
\end{proof}

\bigskip

\begin{rema}\label{rem:modeles}
In our computations made with {\sc Magma}, reflection groups $W$ are represented 
as subgroups of $\Gb\Lb_n(K)$ where $K$ is a number field depending on $W$. 
There are of course infinitely many possibilities for representing $W$ 
in this way, but it turns out that the choice of this model have a considerable 
impact on the time used for computations, and on the form of the defining 
polynomials for the singular varieties we obtain. Let us explain which choices 
we have made and for which reasons:
\begin{itemize}
\item[$\bullet$] We do not use the {\sc Magma} 
command 

\centerline{{\tt ShephardTodd(k)}}

\noindent for defining the complex reflection group $G_k$. 
Indeed, the {\sc Magma} model for $G_k$ is generally not stable under the Galois 
action, and leads to very lengthy computations (and sometimes to computations 
that do not conclude after hours) and to very ugly defining 
polynomials for the singular varieties found by our methods.

\medskip

\item[$\bullet$] In his {\sc Champ} package for {\sc Magma} intended to 
study the representation theory of Cherednik algebras~\cite{thiel}, Thiel 
used the model implemented in the {\sc Chevie} package of {\sc Gap3} 
by Michel~\cite{jean}. These models are almost all stable under the action 
of the Galois group (except for the Coxeter groups $G_{23}=\Wrm(H_3)$ 
and $G_{30}=\Wrm(H_4)$) and leads to much shorter computations 
and much nicer defining polynomials for singular varieties 
(for instance, they almost all have rational coefficients).

\medskip

\item[$\bullet$] We have decided to create our own models for the Coxeter groups 
$G_{23}=\Wrm(H_3)$ and $G_{30}=\Wrm(H_4)$: they are stable under the Galois action 
(so fit with Marin-Michel Theorem). This again shortens the computations 
and lead to polynomials with rational coefficients for defining 
singular varieties: this is how we found en explicit polynomial with rational coefficients 
defining the Sarti dodecic~\cite{bonnafe calculs}. These models are implemented in a file 
{\tt primitive-complex-reflection-groups.m} available in~\cite{bonnafe calculs} 
and are accessible through the command 

\centerline{{\tt PrimitiveComplexReflectionGroup(k)}}

\noindent once this file is downloaded. Note that:
\begin{itemize}
\item[-] This file copies almost entirely Thiel's file except for the Coxeter groups 
$G_{23}=\Wrm(H_3)$, $G_{28}=\Wrm(F_4)$ and $G_{30}=\Wrm(H_4)$.

\item[-] For $G_{23}=\Wrm(H_3)$ and $G_{30}=\Wrm(H_4)$, we have given our own models 
defined over the field $\QM(\r)$, where $\r^4=5\r^2-5$ (i.e. $\r=\sqrt{(5+\sqrt{5})/2}$). 
We do not pretend it is the best possible model but, for our purposes, it is 
the best model available as of today.

\item[-] For $G_{28}=\Wrm(F_4)$, we have used a version which contains 
the Coxeter group $\Wrm(B_4)$ in its standard form (that is, as the group 
of monomial matrices whose non-zero coefficients belong to $\mub_2=\{1,-1\}$) 
as a subgroup of index $3$. This implies in particular that invariant polynomials 
can be expressed in terms on elementary symmetric functions.
\end{itemize}
\end{itemize}
\noindent Of course, as explained in the introduction, the fact that most of the 
singular varieties we construct are defined over $\QM$ do not imply that 
the coordinates of all the singular points are rational, or even real. 
Some of the varieties have in fact no real points. The only example 
where singular points have rational coordinates is given in Example~\ref{ex:f4} 
(see Figure~\ref{fig:f4-8bis}).\finl
\end{rema}

\section{Strategy for finding some ``singular'' invariants in rank $n \ge 3$}
\label{sec:strategy}

\medskip

If $n=2$, then the varieties $\ZC(f)$ are just collections of points, and so are uninteresting 
for our purpose.

\bigskip

\boitegrise{{\bf Hypothesis and notation.} {\it From now on, and until the end of this paper, 
we assume moreover that $n \ge 3$ and that $W$ is {\bfit primitive}. We denote by $r$ the minimal natural 
number such that the space of homogeneous invariants of $W$ of 
degree $d_r$ has dimension $\ge 2$.}}{0.75\textwidth}

\bigskip

Note that this implies that $W$ is one of the groups $G_i$, with $23 \le i \le 37$, 
in Shephard-Todd classification. 
We recall in Table~\ref{tab:degres} the degrees $(d_1,d_2,\dots,d_n)$ of 
these groups. We also give the following informations: the order 
of $W$, the order of $W/\Zrm(W)$ (which is the group which acts faithfully on $\Pb(V)$), 
the degree $d_r$ and, whenever $W$ is a Coxeter group, we recall its type ($\Wrm(X_i)$ denotes the Coxeter 
group of type $X_i$). Recall from general theory that $|W|=d_1d_2\cdots d_n$ and 
$|\Zrm(W)|={\mathrm{Gcd}}(d_1,d_2,\dots,d_n)$. 

\bigskip
\def\petitespace{\vphantom{\DS{\frac{\DS{A}}{\DS{A}}}}}
\def\petitespace{}
\def\espace{}

\begin{table}
$$\begin{array}{!{\vline width 2pt} c !{\vline width 1.2pt} c !{\vline width 1pt} c|c|c|c!{\vline width 2pt}}
\hlinewd{2pt}
\petitespace n & W & |W| & |W/\Zrm(W)| & (d_1,d_2,\dots,d_n) & d_r \\
\hlinewd{2pt}
\multirow{5}{*} 3 &
G_{23}=\Wrm(H_3) & 120 & 60 & 2, 6, 10 & 6 \petitespace\\
\cline{2-6}
& G_{24} & 336 & 168 & 4, 6, 14 & 14 \petitespace\\
\cline{2-6}
& G_{25} & 648 & 108 & 6, 9, 12 & 12 \petitespace\\
\cline{2-6}
& G_{26} & 1~296 & 216 & 6, 12, 18 & 12 \petitespace\\
\cline{2-6}
& G_{27} & 2~160 & 360 & 6, 12, 30 & 12 \petitespace\\
\hlinewd{1.2pt}
\multirow{5}{*} 4 &
G_{28}=\Wrm(F_4) & 1~152 & 576 & 2, 6, 8, 12 & 6 \petitespace\\
\cline{2-6}
& G_{29} & 7~680 & 1~920 & 4, 8, 12, 20 & 8 \petitespace\\
\cline{2-6}
& G_{30}=\Wrm(H_4) & 14~400 & 7~200 & 2, 12, 20, 30 & 12 \petitespace\\
\cline{2-6}
& G_{31} & 46~080 & 11~520 & 8, 12, 20, 24 & 20 \petitespace\\
\cline{2-6}
& G_{32} & 155~520 & 25~920 & 12, 18, 24, 30 & 24 \petitespace\\
\hlinewd{1.2pt}
5 & G_{33} & 51~840 & 25~920 & 4, 6, 10, 12, 18 & 10 \petitespace\\
\hlinewd{1.2pt}
\multirow{2}{*} 6 &
G_{34} & 39~191~040 & 6~531~840 & 6, 12, 18, 24, 30, 42 & 12 \petitespace\\
\cline{2-6}
& G_{35}=\Wrm(E_6) & 51~840 & 25~920 & 2, 5, 6, 10, 12, 14, 18 & 6 \petitespace\\
\hlinewd{1.2pt}
7 & G_{36}=\Wrm(E_7) & 2~903~040 & 1~451~520 & 2, 6, 8, 10, 12, 14, 18 & 6 \petitespace\\
\hlinewd{1.2pt}
8 & G_{37}=\Wrm(E_8) & 696~729~600 & 348~364~800 & 2, 8, 12, 14, 18, 20, 24, 30 & 8 \petitespace\\
\hlinewd{2pt}
\end{array}
$$
\refstepcounter{theo}
\caption{Degrees of primitive complex reflection groups in rank $\ge 3$}\label{tab:degres}
\end{table}

Using {\sc Magma}, we first determine by computer calculations some fundamental invariants 
$f_1$,\dots, $f_r$. By the definition of $r$, the fundamental invariants $f_1$,\dots, $f_{r-1}$ 
are uniquely determined up to scalar. By inspection of Table~\ref{tab:degres}, 
we see that $d_1 < d_2 < \cdots < d_n$ and 
that there is a unique $f$ of the form 
$f_1^{m_1}\cdots f_{r-1}^{m_{r-1}}$ which has degree $d_r$. So the space 
of homogeneous invariants of degree $d_r$ has dimension $2$, and is spanned 
by $f_r$ and $f$. Moreover, all fundamental invariants of degree $d_r$ 
are, up to a scalar, of the form $f_r+u f$, for some $u \in \CM$. 

This means that we need to determine the values of $u$ such that 
$\ZC(f_r+u f)$ is singular. For this, we use the basis $(x_1,\dots,x_n)$ of $V^*$ 
chosen by {\sc Magma} and we set
$$F_u=f_r+u f \qquad\text{and}\qquad F_u^\aff(x_1,\dots,x_{n-1})=F_u(x_1,\dots,x_{n-1},1).$$
This basis allows to identify $\Pb(V)$ with $\Pb^{n-1}(\CM)$ and we denote by 
$\Ab^{n-1}(\CM)$ the affine open subset of $\Pb^{n-1}(\CM)$ defined by $x_n \neq 0$. 
Then $\ZC^\aff(F^\aff_t)$ denotes the affine open subset of $\ZC(F_u)$ defined 
by $x_n \neq 0$. Note the following easy fact:
\equat\label{eq:w-orbit affine}
\text{\it Any $W$-orbit of points in $\Pb^{n-1}(\CM)$ meets $\Ab^{n-1}(\CM)$.}
\endequat
\begin{proof}
Indeed, the linear span of a $W$-orbit of a non-zero vector in $V$ must be equal to $V$, because 
$W$ acts irreducibly. So it cannot be fully contained in the orthogonal of $x_n$.
\end{proof}

\medskip
\def\singfib{{\mathrm{sfib}}}

One deduces immediately the following fact, which will be useful for 
saving much time during computations:
\equat\label{eq:singulier affine}
\text{\it $\ZC(F_u)$ is singular if and only if $\ZC^\aff(F^\aff_u)$ is singular.}
\endequat
Now, let 
$$\XC=\{(\xi,u) \in \Ab^{n-1}(\CM) \times \Ab^1(\CM)~|~F_u^\aff(\xi)=0\}.$$
We denote by $\phi : \XC \to \Ab^1(\CM)$ the second projection. 
Then the fiber $\phi^{-1}(u)$ is the variety $\ZC^\aff(F_u^\aff)$. 
We can then define 
$$\XC_\singfib=\{(\xi,u) \in \XC~|~
\frac{\partial F_u^\aff}{\partial x_1}(\xi) = 
\cdots = \frac{\partial F_u^\aff}{\partial x_{n-1}}(\xi) = 0\}.$$
Then $\XC_\singfib$ is not necessarily the singular locus of $\XC$, but the points 
in $\phi(\XC_\singfib)$ are the values of $u$ for which the fiber $\phi^{-1}(u)=\ZC^\aff(F^\aff_u)$ 
(or, equivalently, $\ZC(F_u)$) 
is singular. We set $U_\sing=\phi(\XC_\singfib)$ and we denote by 
$U_\sing^\irr$ the set of $u \in U_\sing$ such that $\ZC(F_u)$ is irreducible. 
This provides an algorithm for finding these values of $u$: it turns out 
that $\phi$ is not dominant in our examples, so that there are only finitely many such values of $u$. 
We then study more precisely these finite number of cases (number of singular 
points, nature of singularities, Milnor number,\dots). 
Let us see on a simple example how it works:

\bigskip

\begin{exemple}[Coxeter group of type ${\boldsymbol{H_3}}$]\label{ex:h3}
Assume here, and only here, that $W=G_{23}=\Wrm(H_3)$. Then $(d_1,d_2,d_3)=(2,6,10)$ so 
that $r=2$ and $d_r=6$. Then $F_u=f_2 + u f_1^3$. We first define 
$W$ (see Remark~\ref{rem:modeles} for the choice of a model) 
and the fundamental invariants $f_1$ and $f_2$: 

\medskip

\begin{quotation}
\small
\begin{verbatim}
> load 'primitive-complex-reflection-groups.m';
> W:=PrimitiveComplexReflectionGroup(23);
> K<a>:=CoefficientRing(W);
> R:=InvariantRing(W);
> P<x1,x2,x3>:=PolynomialRing(R);
> f1:=InvariantsOfDegree(W,2)[1];
> f2:=InvariantsOfDegree(W,6)[1];
> Gcd(f1,f2);
1
\end{verbatim}
\end{quotation}

\medskip

\noindent Note that the last command shows that the invariant $f_2$ of degree $6$ 
we have chosen is indeed a fundamental invariant. We now define $F_u^\aff$ 
and $\XC_\singfib$ and then 
determine the set $U_\sing$ of values of $u$ such that $\ZC(F_u)$ is singular:

\medskip

\begin{quotation}
\small
\begin{verbatim}
> P2:=Proj(P);
> A2xA1<xx1,xx2,u>:=AffineSpace(K,3);
> A1<U>:=AffineSpace(K,1);
> phi:=map<A2xA1->A1 | [u]>;
> f1aff:=Evaluate(f1,[xx1,xx2,1]);
> f2aff:=Evaluate(f2,[xx1,xx2,1]);
> Fuaff:=f2aff + u * f1aff^3;
> X:=Scheme(A2xA1,Fuaff);
> Xsfib:=Scheme(X,[Derivative(Fuaff,i) : i in [1,2]]);
> Psing:=MinimalBasis(phi(Xsfib));
> # Psing;
1
> Factorization(Psing[1]);
[
    <T + 1, 1>,
    <T + 9/10, 1>,
    <T + 63/64, 1>
]
> Using:=[-1, -9/10, -63/64];
\end{verbatim}
\end{quotation}

\medskip

\noindent We next determine for which values $u \in U_\sing=\{u_1,u_2,u_3\}$ 
the curve $\ZC(F_u)$ is irreducible:

\medskip

\begin{quotation}
\small
\begin{verbatim}
> F:=[f2+ui*f1^3 : ui in Using]; // the polynomials F_{t_i}
> Z:=[Curve(P2,f) : f in F];
> [IsAbsolutelyIrreducible(i) : i in Z];
[ true, true, false ]
\end{verbatim}
\end{quotation}

\noindent We then study the singular locus of the irreducible curves 
$\ZC(F_u)$ for $u=u_1$ or $u_2$. Let us see how to do it for $u=u_1$:

\begin{quotation}
\small
\begin{verbatim}
> Z1sing:=SingularSubscheme(Z[1]);   
> Z1sing:=ReducedSubscheme(Z1sing);
> Degree(Z1sing);
10
> points:=SingularPoints(Z[1]);
> # points;
10
> pt:=points[1];
> IsNode(Z[1],pt);
true
> # ProjectiveOrbit(W,pt);
10
\end{verbatim}
\end{quotation}

\medskip

\noindent The command {\tt Degree(Z1sing)} 
shows that $\ZC(F_{u_1})$ contains exactly $10$ singular points. 
The command {\tt \# points} shows that they are all defined over 
the field {\tt K} ($=\QM(\sqrt{5})$). The command {\tt \# ProjectiveOrbit(W,p1)} shows that they 
are all in the same $W$-orbit (the function {\tt ProjectiveOrbit} 
has been defined by the author for computing orbits in projective 
spaces (see~\cite{1440} or~\cite{bonnafe calculs} for the code). So all these singularities 
are equivalent and the command {\tt IsNode(Z[1],pt)} shows that 
they are all nodes.

One can check similarly that $\ZC(F_{u_2})$ has $6$ nodes, all belonging 
to the same $W$-orbit.\finl
\end{exemple}

\bigskip

In the next sections, we will give tables of singular curves and surfaces obtained in this way. 
Inspection of these tables (and Examples~\ref{ex:g33} and~\ref{ex:e6}) leads to the following result:

\bigskip

\begin{prop}\label{prop:q-all}
Apart from the two singular surfaces $\SC$ and $\SC'$ of degree $8$ with $80$ nodes defined 
by invariants of $G_{29}$, all the singular curves and surfaces described in 
Tables~\ref{tab:singular curves},~\ref{tab:g26} and~\ref{tab:singular surfaces} can be defined over $\QM$. 
The singular surfaces $\SC$ and $\SC'$ are Galois conjugate over $\QM$.
\end{prop}

\bigskip

\begin{proof}
One could just check that the polynomials given thanks to the {\sc Magma} codes 
contained in~\cite{bonnafe calculs} have coefficients in $\QM$. But 
one could also follow the same argument as in Proposition~\ref{prop:q}, based on Marin-Michel Theorem, 
by using the fact that all these singular curves and surfaces are characterized 
by their number of singular points or their type. 
\end{proof}

\bigskip

\begin{prop}\label{prop:transitive}
If $W=G_k$, with $23 \le k \le 35$ and $k \neq 34$, and if $u \in U_\sing^\irr$, 
then $W$ acts transitively on $\ZC_\sing(F_u)$.
\end{prop}

\bigskip

\section{Singular curves from groups of rank $3$}

\medskip

\boitegrise{{\bf Hypothesis.} {\it We still assume that $W$ is primitive but, in this section, 
we assume moreover that $n=3$.}}{0.75\textwidth}

\medskip

This means that $W$ is one of the groups $G_i$, for $23 \le i \le 27$. 
We denote by $(f_1,f_2,f_3)$ a set of fundamental invariants provided by {\sc Magma}. 
Table~\ref{tab:singular curves} gives the list of curves obtained through the methods detailed 
in Section~\ref{sec:strategy}. This table contains the 
degree $d_r$, the cardinality of $U_\sing^\irr$, 
the number of singular points and the type of the singularity 
(since all singular points belong to the same $W$-orbit by Proposition~\ref{prop:transitive}, 
they are all equivalent singularities). 
Details of {\sc Magma} computations are given in~\cite{bonnafe calculs} 
(they follow the lines of Example~\ref{ex:h3}). We use standard notation for the types of the 
singularities of curves~\cite{agv}. 
For instance (here, we denote by $m$ the multiplicity, $\mu$ the Milnor number and 
$\t$ the Tjurina number):
\begin{itemize}
\item[$\bullet$] $A_1$ is a {\it node}, i.e. a singularity equivalent to $xy$: 
in this case, $m=2$ and $\mu=\t=1$.

\item[$\bullet$] $A_2$ is a {\it cusp}, i.e. a singularity equivalent to $y^2-x^3$: 
in this case, $m=2=\mu=\t=2$.

\item[$\bullet$] $D_4$ is a singularity equivalent to $x(y^2-x^2)$: in this 
case, $m=3$, $\mu=\t=4$.

\item[$\bullet$] $X_9$ is a singularity equivalent to $xy(y-x)(y+x)$: 
in this case, $m=4$, $\mu=\t=9$.

\item[$\bullet$] $E_6$ is a singularity equivalent to $y^3-x^4$: in this case, 
$m=3$, $\mu=\t=6$.
\end{itemize}

\bigskip

\bigskip

\def\espace{\vphantom{\DS{\frac{\DS{A^A}}{\DS{A_A}}}}}

\begin{table}
$$\begin{array}{!{\vline width 2pt} c !{\vline width 1.2pt} c !{\vline width 1pt} c|c|c|c!{\vline width 2pt}}
\hlinewd{2pt}
\petitespace W & d_r & |U_\sing^\irr| & u_i & |\ZC_\sing(F_{u_i})| & \text{Singularity} \\
\hlinewd{2pt}
\multirow{2}{*}{$G_{23}=\Wrm(H_3)\espace$} &\multirow{2}{*}{$6\espace$} & 
\multirow{2}{*}{$2\espace$}
& u_1 & 6 & A_1 \petitespace\\
\cline{4-6}
&&& u_2 & 10 & A_1  \petitespace\\
\hlinewd{1.2pt}
\multirow{3}{*}{$G_{24}\espace$} &\multirow{3}{*}{$14\espace$} & \multirow{3}{*}{$3 \espace$}
& u_1 & 21 & A_1 \petitespace\\
\cline{4-6}
&&& u_2 & 28 & A_1  \petitespace\\
\cline{4-6}
&&& u_3 & 42 & A_2  \petitespace\\
\hlinewd{1.2pt}
\multirow{2}{*}{$G_{25}\espace$} &\multirow{2}{*}{$12\espace$} & \multirow{2}{*}{$2 \espace$}
& u_1 & 12 & D_4\petitespace\\
\cline{4-6}
&&& u_2 & 36 & A_2  \petitespace\\
\hlinewd{1.2pt}
\multirow{2}{*}{$G_{27}\espace$} &\multirow{2}{*}{$12\espace$} & \multirow{2}{*}{$2 \espace$}
& u_1 & 45 & A_1 \petitespace\\
\cline{4-6}
&&& u_2 & 36 & A_1  \petitespace\\
\hlinewd{2pt}
\end{array}
$$
\refstepcounter{theo}
\caption{Singularities of the curves $\ZC(F_u)$ for $t \in U_\sing^\irr$}\label{tab:singular curves}
\end{table}

\begin{exemple}\label{ex:record?}
A plane curve is called {\it cuspidal} if all its singular points are of type $A_2$. 
By~\cite[(0.4)]{sakai}, a cuspidal plane curve of degree $14$ has at most $55$ singular points of type $A_2$. 
But it is not known if this is the sharpest bound: to the best of our knowledge, 
no cuspidal plane curve of degree $14$ 
with $42$ or more singular points of type $A_2$ was known before the above example 
of $\ZC(F_{u_3})$ for $W=G_{24}$.

Also, a cuspidal plane curve of degree $12$ can have at most $40$ singular points~\cite[(0.4)]{sakai}, 
but it is not known if this bound can be achieved. However, there exists at least one cuspidal 
curve of degree $12$ with $39$ cusps~\cite[Example~6.3]{39}. Our example obtained 
from invariants of $G_{25}$, with $36$ cusps, approaches these bounds and has an automorphism 
group of order $\ge 108$.\finl
\end{exemple}

\bigskip

\begin{rema}\label{rem:g26}
Note that $G_{26}$ does not appear in Table~\ref{tab:singular curves}. The reason is the following: 
if $W=G_{26}$, then $d_r=12$ but $G_{26}$ contains $W'=G_{25}$ as a normal subgroup of index $2$ 
and it turns out that invariants of degree $12$ of $G_{25}$ and $G_{26}$ coincide. 
This makes the computation for $G_{26}$ unnecessary in this case. Note, however, 
the next Example~\ref{ex:g26}, where we construct singular curves 
of degree $18$ using invariants of $G_{26}$.\finl
\end{rema}

\bigskip

\begin{exemple}[The group ${\boldsymbol{G_{26}}}$]\label{ex:g26}
We assume in this example that $W=G_{26}$. Recall that $(d_1,d_2,d_3)=(6,12,18)$. 
Up to a scalar, any fundamental invariant of degree $18$ of $W$ is of the form $F_{u,v}=f_3+uf_1f_2+vf_1^3$ 
for some $(u,v) \in \Ab^2(\CM)$. 
Using {\sc Magma}, one can check the following facts. First, the set $\CC$ of 
$(u,v) \in \Ab^2(\CM)$ such that 
$\ZC(F_{u,v})$ is singular is a union of three affine lines $\LC_1$, $\LC_2$, $\LC_3$ 
and a smooth curve $\EC$ isomorphic to $\Ab^1(\CM)$. The singular locus $\CC_\sing$ 
of $\CC$ consists of $7$ points and it turns out that there are only $5$ points 
$(u_i,v_i)_{1 \le i \le 5}$ in $\CC_\sing$ such that $\ZC(F_{u_i,v_i})$ is 
irreducible. Table~\ref{tab:g26} gives the information about singularities 
of these varieties $\ZC(F_{u_i,v_i})$ (with the numbering used 
in our {\sc Magma} programs~\cite{bonnafe calculs}).

\begin{table}
$$\begin{array}{!{\vline width 2pt} c !{\vline width 1.2pt} c !{\vline width 1pt} c|c!{\vline width 2pt}}
\hlinewd{2pt}
(u,v) & |\ZC_\sing(F_{u,v})| & W{\operatorname{\!-}}\text{orbits} & \text{Singularity} \petitespace\\
\hlinewd{2pt}
\multirow{2}{*}{$(u_1,v_1)\espace$} & \multirow{2}{*}{$63\espace$} & 9 & 
X_9 \petitespace\\
\cline{3-4} && 54 & A_2 \\
\hlinewd{1.2pt}
\multirow{2}{*}{$(u_2,v_2)\espace$} & \multirow{2}{*}{$21\espace$} & 9 & 
X_9 \petitespace\\
\cline{3-4} && 12 & D_4 \petitespace \\
\hlinewd{1.2pt}
\multirow{2}{*}{$(u_3,v_3)\espace$} & \multirow{2}{*}{$45\espace$} & 9 & 
X_9 \petitespace\\
\cline{3-4}
&& 36 & A_2 \petitespace \\
\hlinewd{1.2pt}
(u_4,v_4) & 36 & 36 & E_6 \petitespace\\
\hlinewd{1.2pt}
\multirow{2}{*}{$(u_5,v_5)\espace$} & \multirow{2}{*}{$84\espace$} & 12 & 
D_4 \petitespace\\
\cline{3-4}
&& 72 & A_2 \petitespace \\
\hlinewd{2pt}
\end{array}$$
\refstepcounter{theo}
\caption{Some singular curves of degree $18$ defined by invariants of $G_{26}$}\label{tab:g26}
\end{table}
Note that a cuspidal curve of degree $18$ 
has at most $94$ singularities of type $A_2$~\cite[(0.3)]{sakai}. Note also that there exists a cuspidal 
curve of degree $18$ with $81$ cusps~\cite{ivinskis}.\finl
\end{exemple}

\bigskip

\section{Singular surfaces from groups of rank $4$}\label{sec:rang 4}

\medskip

\boitegrise{{\bf Hypothesis.} {\it We still assume that $W$ is primitive but, in this section, 
we assume moreover that $n=4$.}}{0.75\textwidth}

\medskip

This means that $W$ is one of the groups $G_i$, for $28 \le i \le 32$. 
We denote by $(f_1,f_2,f_3,f_4)$ a set of fundamental invariants provided by {\sc Magma} 
and we denote by $U_\sing^\irr$ the set of elements $u \in \CM$ such that 
$\ZC(F_u)$ is irreducible and singular. 
Table~\ref{tab:singular surfaces} gives the list of surfaces obtained through the methods detailed 
in Section~\ref{sec:strategy}. 
This table contains the degree $d_r$, the 
number of values of $t$ such that $\ZC(F_u)$ is irreducible and singular, 
the number of singular points and informations about the singularity 
(since all singular points belong to the same $W$-orbit by Proposition~\ref{prop:transitive}, 
they are all equivalent singularities). The number $m$ (resp. $\mu$, resp. $\t$) denotes 
the multiplicity (resp. the Milnor number, resp. the Tjurina number). 

\begin{table}
$$\begin{array}{!{\vline width 2pt} c !{\vline width 1.2pt} c !{\vline width 1pt} c|c|c|c!{\vline width 2pt}}
\hlinewd{2pt}
\petitespace W & d_r & |U_\sing^\irr| & u_i & |\ZC_\sing(F_{u_i})| & \text{Singularity} \\
\hlinewd{2pt}
\multirow{4}{*}{$G_{28}=\Wrm(F_4)\espace$} &\multirow{4}{*}{$6\espace$} & 
\multirow{4}{*}{$4\espace$}
& u_1 & 12 & A_1 \petitespace\\
\cline{4-6}
&&& u_2 & 12 & A_1  \petitespace\\
\cline{4-6}
&&& u_3 & 48 & A_1  \petitespace\\
\cline{4-6}
&&& u_4 & 48 & A_1  \petitespace\\
\hlinewd{1.2pt}
\multirow{5}{*}{$G_{29}\espace$} &\multirow{5}{*}{$8\espace$} & \multirow{5}{*}{$5 \espace$}
& u_1 & 40 & A_1 \petitespace\\ 
\cline{4-6}
&&& u_2 & 20 & \text{Ordinary}, m=3, \mu=11, \tau=10 \petitespace\\
\cline{4-6}
&&& u_3 & 160 & A_1  \petitespace\\
\cline{4-6}
&&& u_4 & 80 & A_1  \petitespace\\
\cline{4-6}
&&& u_5 & 80 & A_1  \petitespace\\
\hlinewd{1.2pt}
\multirow{4}{*}{$G_{30}=\Wrm(H_4)\espace$} &\multirow{4}{*}{$12\espace$} & \multirow{4}{*}{$4 \espace$}
& u_1 & 300 & A_1\petitespace\\
\cline{4-6}
&&& u_2 & 60 & A_1  \petitespace\\
\cline{4-6}
&&& u_3 & 360 & A_1  \petitespace\\
\cline{4-6}
&&& u_4 & 600 & A_1  \petitespace\\
\hlinewd{1.2pt}
\multirow{5}{*}{$G_{31}\espace$} &\multirow{5}{*}{$20\espace$} & \multirow{5}{*}{$5 \espace$}
& u_1 & 480 & A_1 \petitespace\\
\cline{4-6}
&&& u_2 & 960 & A_1  \petitespace\\
\cline{4-6}
&&& u_3 & 1~\!920 & A_1  \petitespace\\
\cline{4-6}
&&& u_4 & 640 & A_1  \petitespace\\
\cline{4-6}
&&& u_5 & 1~\!440 & A_1  \petitespace\\
\hlinewd{1.2pt}
\multirow{4}{*}{$G_{32}\espace$} &\multirow{4}{*}{$24\espace$} 
& \multirow{4}{*}{$4 \espace$} & u_1 & 40 &  \text{Ordinary}, m=6, \mu=125, \t=125\petitespace\\
\cline{4-6}
&&& u_2 & 360 & \text{Non-ordinary}, m=3, \mu=18, \t=18  \petitespace\\
\cline{4-6}
&&& u_3 & 1~\!440 & D_4 \petitespace\\
\cline{4-6}
&&& u_4 & 540 & \text{Non-simple, non-ordinary}, m=2, \mu=9, \t=9  \petitespace\\
\hlinewd{2pt}
\end{array}
$$
\refstepcounter{theo}
\caption{Singularities of the surfaces $\ZC(F_u)$ for $u \in U_\sing^\irr$}\label{tab:singular surfaces}
\end{table}

\bigskip

The example with $1~\!440$ singularities of type $D_4$ 
obtained from $G_{32}$ is detailed in section~\ref{sec:g32}: one can derive from 
the construction a surface of degree $8$ with $44$ singularities of type $D_4$ 
(see also~\cite{1440}). 

\bigskip

\begin{rema}[Coxeter groups of rank ${\boldsymbol{4}}$]\label{rem:f4-h4}
In Table~\ref{tab:singular surfaces}, 
the cases of Coxeter groups of type $F_4$ and $H_4$ (i.e. the primitive 
reflection groups $G_{28}$ and $G_{30}$) was dealt with by Sarti~\cite{sarti 0}.\finl
\end{rema}

\bigskip

\begin{center}
\begin{figure}
\begin{center}
\includegraphics[scale=0.2]{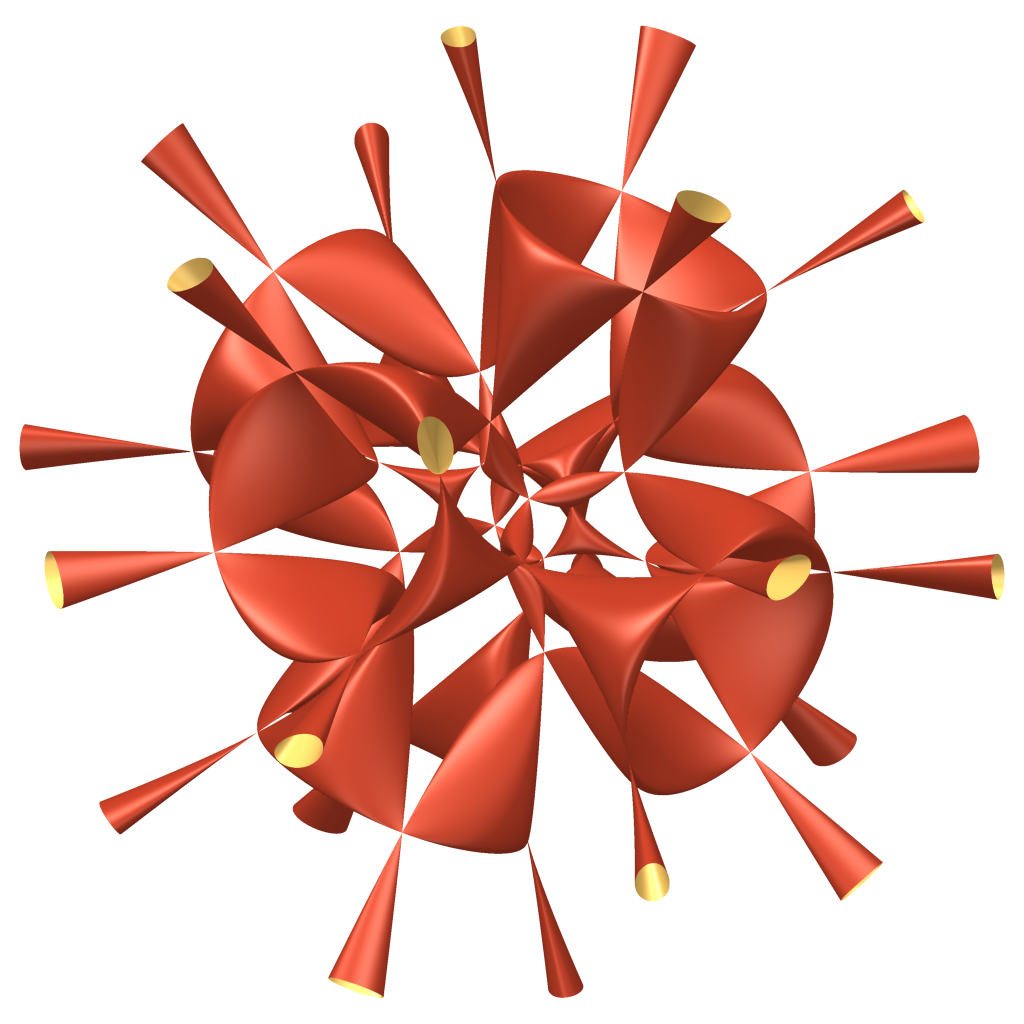}
\includegraphics[scale=0.2]{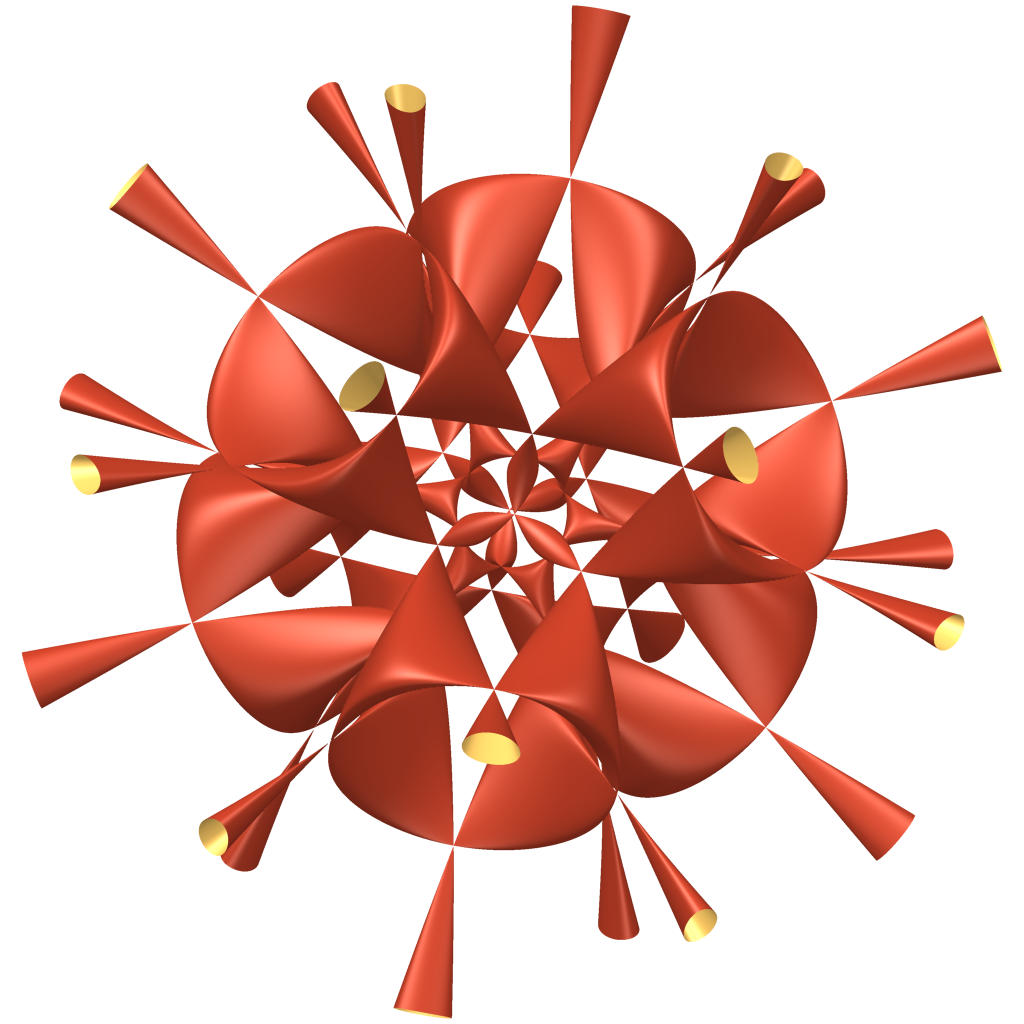}
\caption{Part of the real locus of $\ZC(\ph_2)$ for $W=G_{28}=\Wrm(F_4)$.}\label{fig:f4-8bis}
\end{center}
\end{figure}
\end{center}

\begin{exemples}[Coxeter group of type ${\boldsymbol{F_4}}$]\label{ex:f4}
Assume in this example, and only in this example, that $W=G_{28}=\Wrm(F_4)$ is the Coxeter group 
of type $F_4$, in the form explained in Remark~\ref{rem:modeles}. 
We denote by $\s_1$, $\s_2$, $\s_3$, $\s_4$ the elementary symmetric polynomials 
in $x_1$, $x_2$, $x_3$, $x_4$ and if $f \in \CM[x_1,x_2,x_3,x_4]$ and $k \ge 1$, 
we set $f[k]=f(x_1^k,x_2^k,x_3^k,x_4^k)$. 

Let $\ph_1$ and $\ph_2$ be the following two polynomials:
$$\ph_1=7\sigma_1[2]^4-72\sigma_1[2]^2\sigma_2[2]+4~\!320\sigma_4[2]+432\s_2[4]$$
$$\ph_2=\sigma_1[2]^4-9\sigma_1[2]^2\s_2[2]+27\s_2[2]^2-27\s_1[2]\s_3[2]+324\s_4[2].\leqno{\text{and}}$$
Then it is easily checked that $\ph_i \in \CM[V]^W$ 
and that the two varieties $\ZC(\ph_i)$ are isomorphic (because there is an element $g$ of 
$\Nrm_{\Gb\Lb_4(\CM)}(W)$ such that 
$\ph_2=\lexp{g}{\ph_1}$) and have the following properties:
\begin{itemize}
\item[$\bullet$] The reduced singular locus 
$\ZC_\sing(\ph_i)$ has dimension $0$ and consists of $48$ points 
which are all quotient singularities of type $D_4$.

\item[$\bullet$] The group $G_{28}$ acts transitively on $\ZC_\sing(\ph_i)$ and all elements 
of $\ZC_\sing(\ph_i)$ have coordinates in $\QM$.
\end{itemize}
This shows in particular that 
\equat\label{eq:mud48}
\mu_{D_4}(8) \ge 48,
\endequat
as announced in the introduction. 
Figure~\ref{fig:f4-8bis} shows part of the real locus of $\ZC(\ph_2)$.\finl
\end{exemples}

\bigskip

%
%

\bigskip

\begin{center}
\begin{figure}
\begin{center}
\includegraphics[scale=0.2]{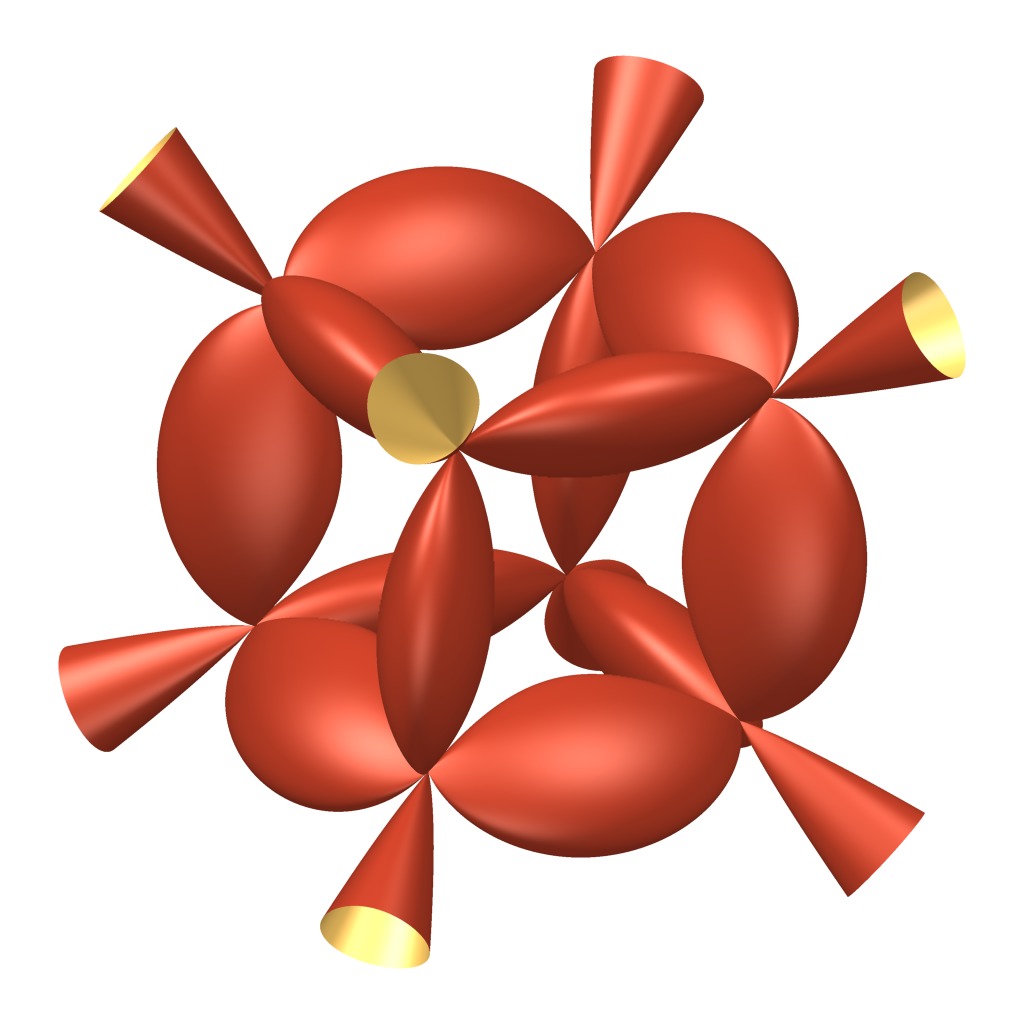}
\includegraphics[scale=0.2]{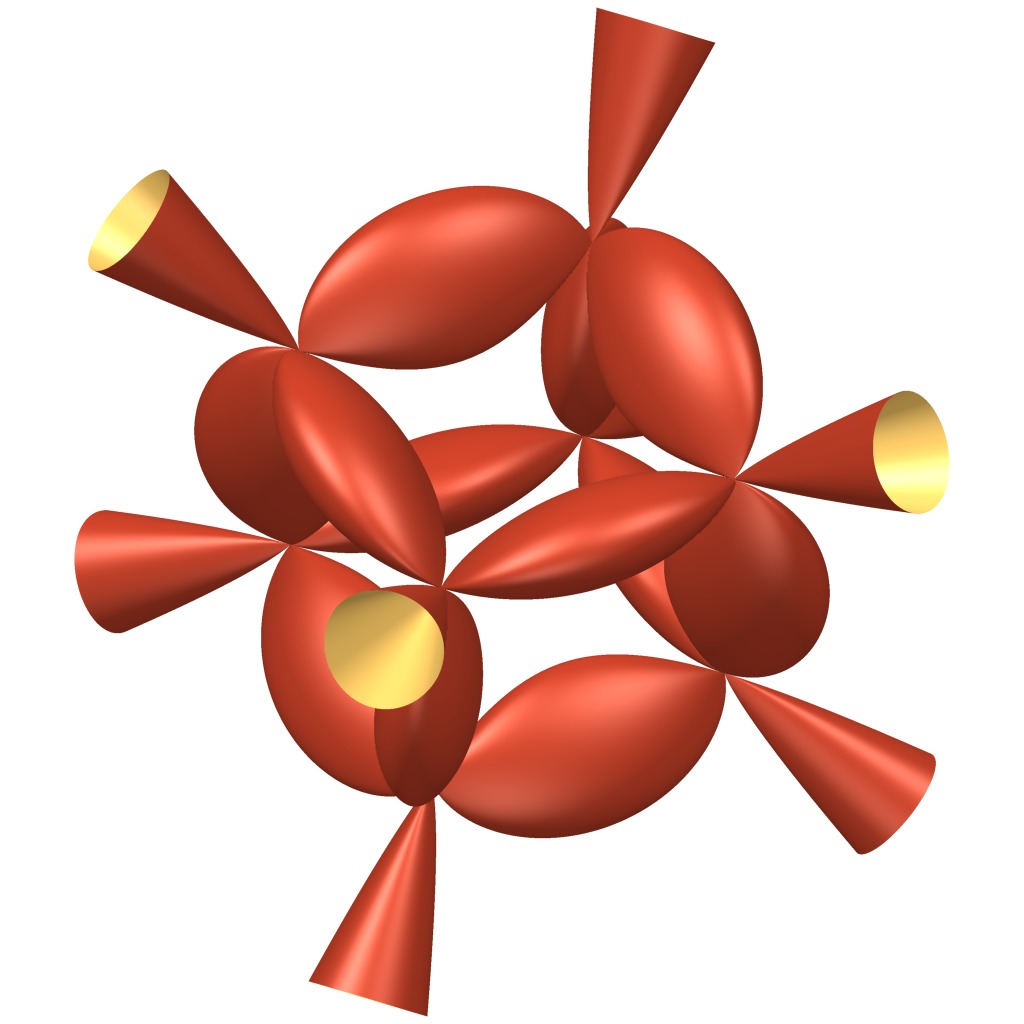}
\caption{Part of the real locus of $\ZC(F_{u_3})$ for $W=G_{29}$.}\label{fig:g29}
\end{center}
\end{figure}
\end{center}

\begin{center}
\begin{figure}
\begin{center}
\includegraphics[scale=0.2]{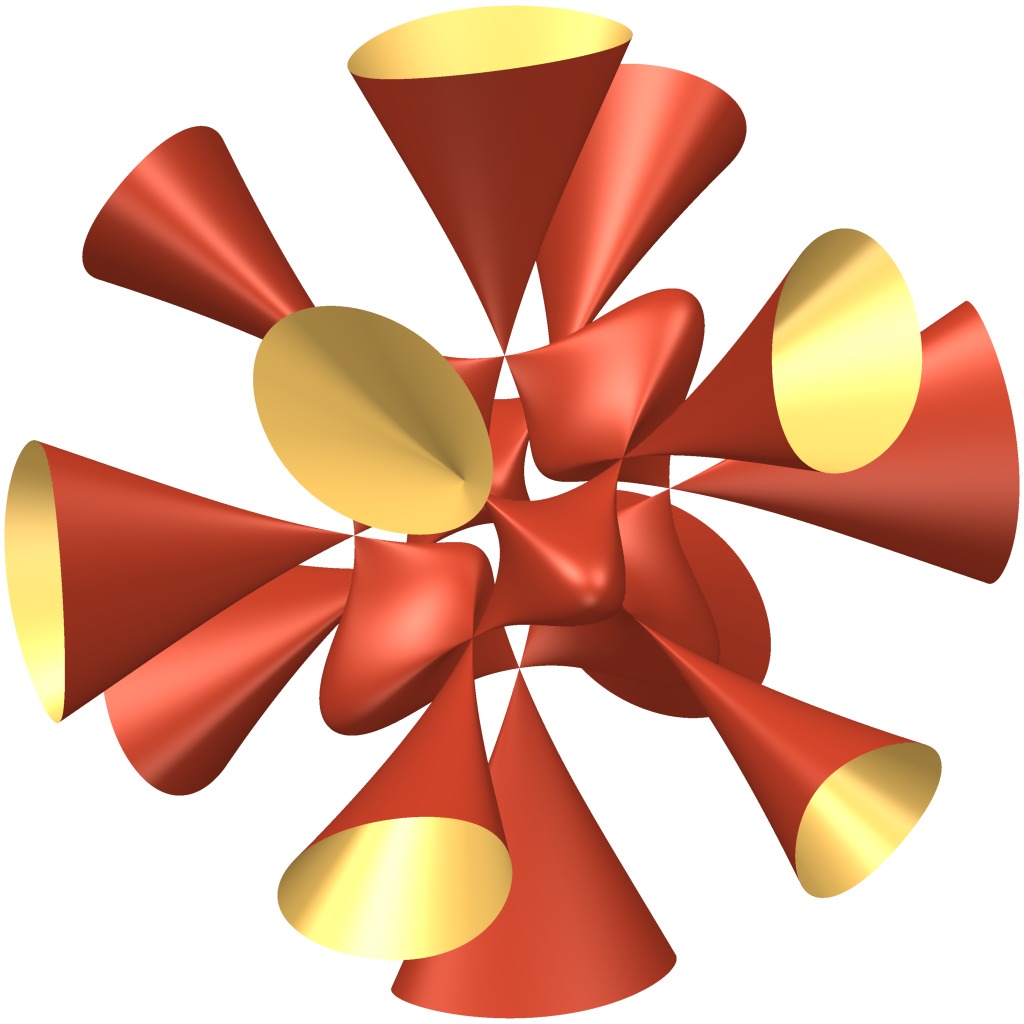}
\includegraphics[scale=0.2]{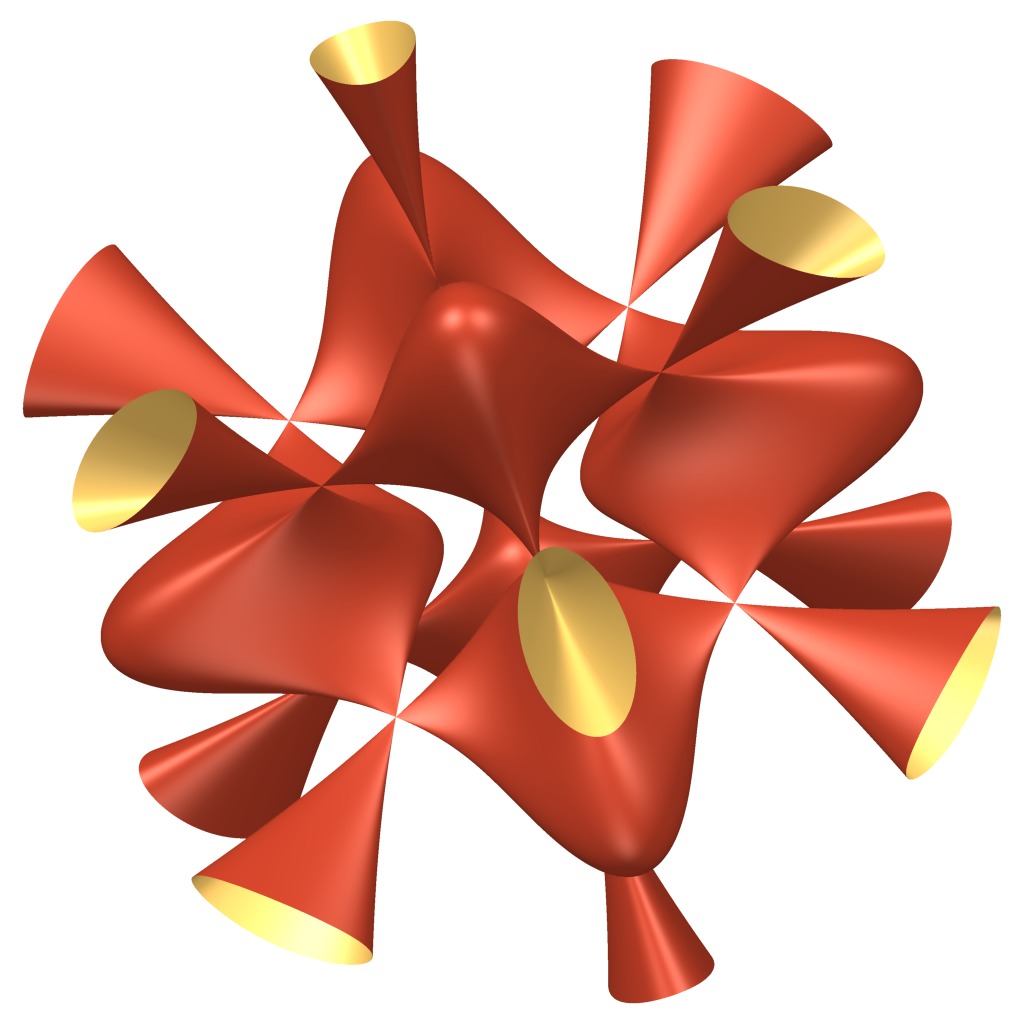}
\caption{Part of the real locus of $\ZC(\ph_1)$ for $W=G_{29}$.}\label{fig:g29-12}
\end{center}
\end{figure}
\end{center}

\begin{center}
\begin{figure}
\begin{center}
\includegraphics[scale=0.2]{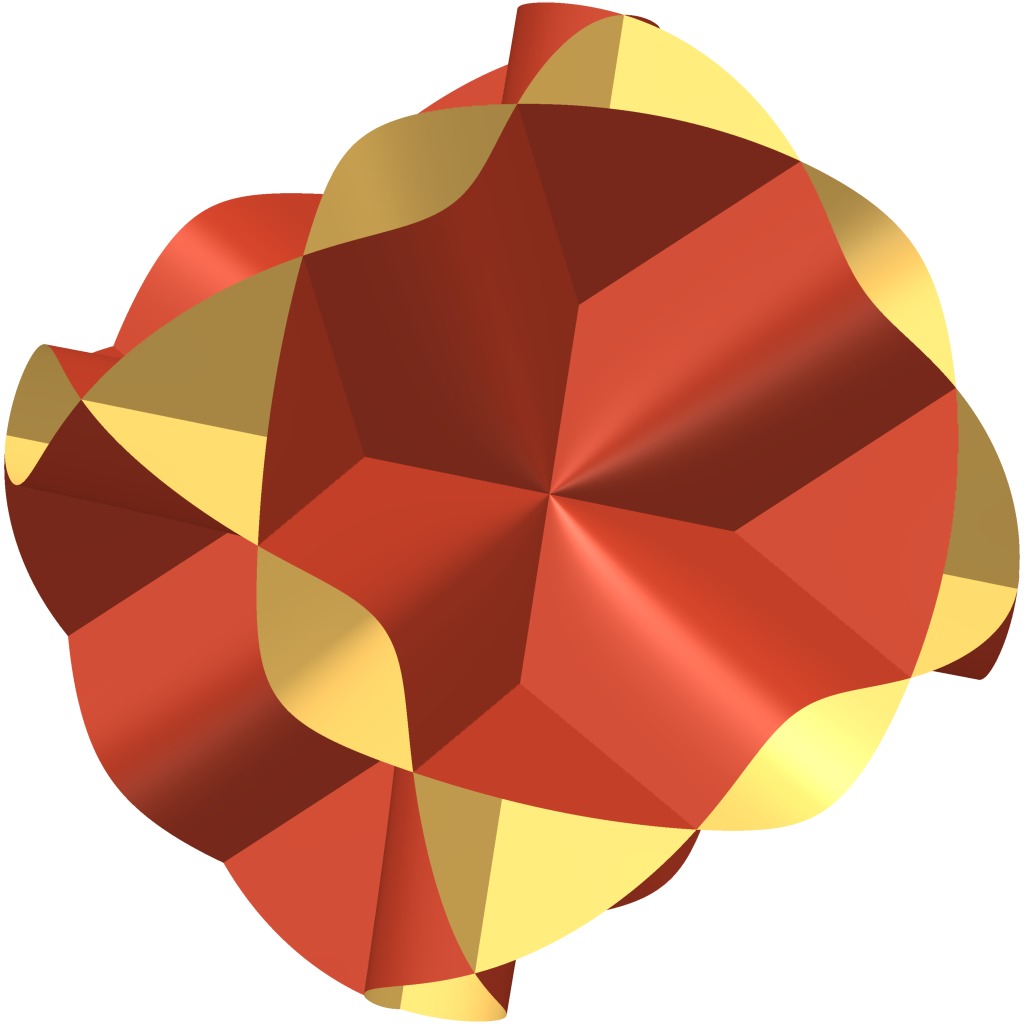}
\includegraphics[scale=0.2]{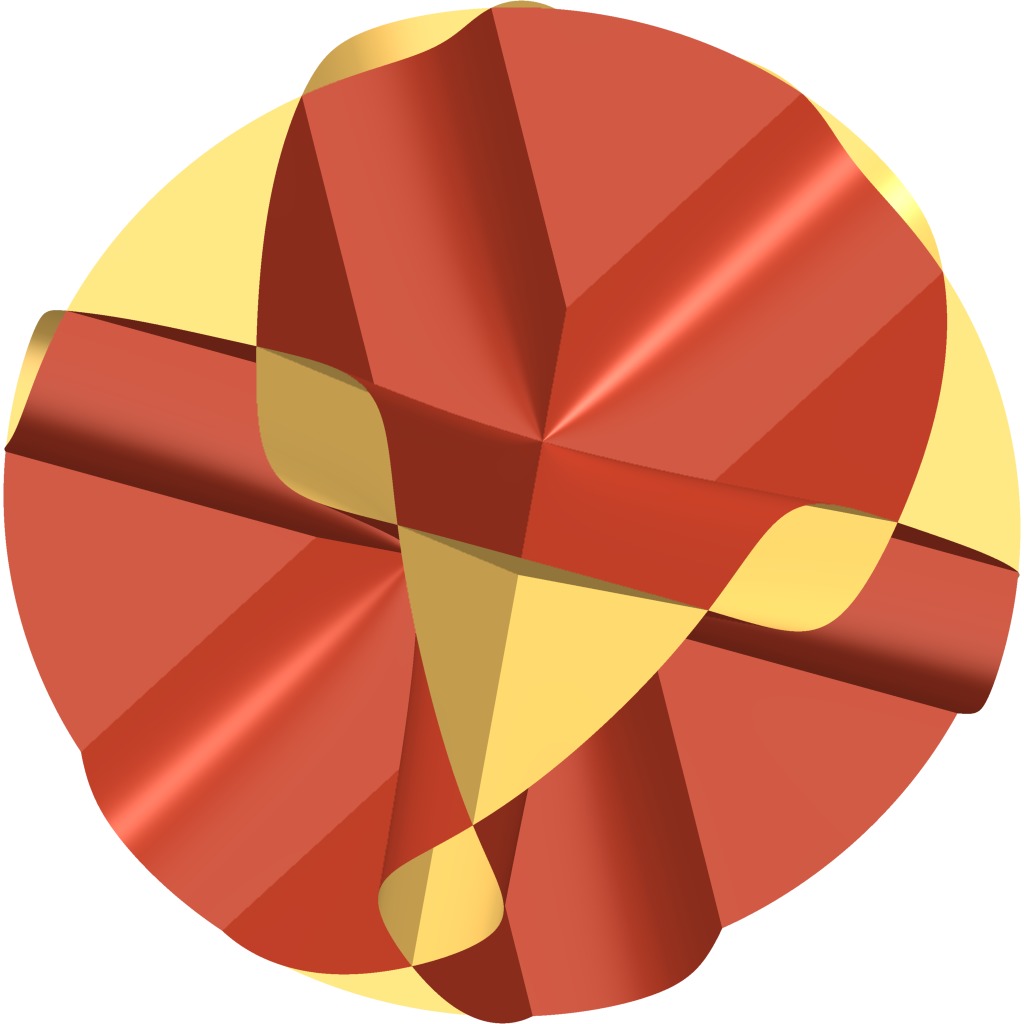}
\caption{Part of the real locus of $\ZC(\ph_2)$ for $W=G_{29}$.}\label{fig:g29-zsing1}
\end{center}
\end{figure}
\end{center}

\begin{exemples}[The group ${\boldsymbol{G_{29}}}$]\label{ex:g29}
Assume in this example, and only in this example, that $W=G_{29}$, in the version implemented 
by Jean Michel in the {\tt Chevie} package of {\tt GAP3}~\cite{jean}. 
Then it contains the symmetric group $\SG_4$ (viewed as the subgroup 
of $\Gb\Lb_4(\CM)$ consisting of permutation matrices). We use the 
notation of Example~\ref{ex:f4} for elementary symmetric functions 
and evaluation at powers of the indeterminates. 

\medskip

(1) Recall that the {\it Endra\ss~octic}~\cite{endrass}  
has degree $8$ and $168$ nodes and its automorphism group 
has order $16$. As shown in Table~\ref{tab:singular surfaces}, 
$\ZC(F_{u_3})$ is an irreducible surface in $\Pb^3(\CM)$ 
with $160$ nodes and a group of automorphisms of order at least $1~\!920$, 
thus approaching Endra\ss'~record but with more symmetries. However, this surface has no real point. 
Up to a scalar, we have 
$$F_{u_3}=\sigma_1[8]+3\s_1[2]^2\s_2[2]+2\s_2[4]-30\s_1[2]\s_3[2]+240\s_4[2].$$
It is still an open 
question to determine whether one can find a surface of degree $8$ 
in $\Pb^3(\CM)$ with more than $168$ nodes (being aware that the 
maximal number of nodes cannot exceed $174$, see~\cite{miyaoka}). 

\medskip

(2) For the surface $\ZC(F_{u_2})$, it can be shown with the software {\sc Singular} 
that the singularities are all of type $T_{4,4,4}$ that is, are equivalent 
to the singularity $xyz+x^4+y^4+z^4$. Up to a scalar, we have
$$F_{u_2}=\s_1[2]^4 -32 \s_1[2]\s_3[2]+256\s_4[2].$$
Figure~\ref{fig:g29} shows part of the real locus of $\ZC(F_{u_2})$.

\medskip

(3) On the other hand, if we set 
$$\ph_1=
\sigma_1[2]^6-\frac{3}{2}\sigma_1[2]^4\sigma_2[2]-78\sigma_1[2]^2\sigma_2[2]^2
+\frac{585}{2}\sigma_1[2]^3\sigma_3[2]+208\sigma_2[2]^3$$
$$-990\sigma_1[2]\sigma_2[2]\sigma_3[2]
+1710\sigma_1[2]^2\sigma_4+1350\sigma_3[2]^2-2880\sigma_2[2]\sigma_4[2],$$
we can check that $\ph_1 \in \CM[V]^W$ and that:
\begin{itemize}
\item[$\bullet$] $\ZC(\ph_1)$ has exactly $160$ singular points, 
which are all singularities of type $D_4$.

\item[$\bullet$] $\ZC_\sing(\ph_1)$ is a single $G_{29}$-orbit.
\end{itemize}
This shows that
\equat
\mu_{D_4}(12) \ge 160,
\endequat
as announced in the introduction. 
This improves considerably known lower bounds (to the best of our knowledge, it 
was only known that $\mu_{D_4}(12) \ge 96$, see~\cite{escudero 2}). 
Recall also that Miyaoka's bound says that $\mu_{D_4}(12) \le 198$. Figure~\ref{fig:g29-12} 
shows part of the real locus of $\ZC(\ph_1)$.

\medskip

(4) Let us keep going on with fundamental invariants of degree $12$. 
Let 
$$\ph_2=\sigma_3[2]\sigma_1[2]^3-4\sigma_1[2]\sigma_2[2]\sigma_3[2]
+4\sigma_1[2]^2\sigma_4[2]+4\sigma_3[2]^2$$
(up to a scalar). Then 
$\ph_2 \in \CM[V]^W$ is irreducible over $\CM$ (this has been checked with {\sc Singular}) 
and computations with {\sc Magma} show that:
\begin{itemize}
\item[$\bullet$] $\ZC_\sing(\ph_2)$ has pure dimension $1$ and is the 
union of $30$ lines. 

\item[$\bullet$] $G_{29}$ acts transitively on these $30$ lines.

\item[$\bullet$] The set of points belonging to at least two of these $30$ lines 
has cardinality $60$, and splits into two $G_{29}$-orbits (one of 
cardinality $40$, the other of cardinality $20$). 
\end{itemize}
Figure~\ref{fig:g29-zsing1} shows part of the real locus of $\ZC(\ph_2)$.\finl
\end{exemples}

\bigskip

\bigskip

\bigskip

\bigskip

\begin{exemple}[The group ${\boldsymbol{G_{31}}}$]\label{ex:g31}
Recall that the {\it Chmutov surface}~\cite{chmutov} of degree $20$ has $2~\!926$ nodes and 
that an irreducible surface in $\Pb^3(\CM)$ of degree $20$ cannot have more than $3~\!208$ nodes~\cite{miyaoka}. 
The third surface associated with $G_{31}$ in Table~\ref{tab:singular surfaces} 
has ``only'' $1~\!920$ nodes and most of them are not real (contrary to the Chmutov surface). 
However, it has a big group of automorphisms (of order a least $11~\!520$).\finl
\end{exemple}

\bigskip

\section{Examples in higher dimension}

\bigskip

\begin{exemple}[The group ${\boldsymbol{G_{33}}}$]\label{ex:g33}
Computations with {\sc Magma} show that there are no fundamental invariant $f_3$ 
of degree $10$ of $G_{33}$ such that $\ZC(f_3)$ is singular~\cite{bonnafe calculs}.\finl
\end{exemple}

\bigskip

\begin{exemple}[Coxeter group of type ${\boldsymbol{E_6}}$]\label{ex:e6}
Assume in this Example, and only in this Example, that $W=G_{35}$ is a Coxeter group 
of type $E_6$. Then $r=3$ and $(d_1,d_2,d_3)=(2,5,6)$, so that 
any fundamental invariant of degree $6$ of $W$ is of the form $F_u=f_3+u f_1^3$ 
for some $u \in \CM$. Computations with {\sc Magma} show that~\cite{bonnafe calculs}:
\begin{itemize}
\itemth{a} $U_\sing=U_\sing^\irr$ has cardinality $8$. 

\itemth{b} For each $u \in U_\sing$, $\ZC_\sing(F_u)$ has dimension $0$, 
$W$ acts transitively on $\ZC_\sing(F_u)$, and all these singular points 
are nodes.

\itemth{c} The hypersurfaces $\ZC(F_u)$, $u \in U_\sing^\irr$, have respectively 
$27$, $36$, $135$, $216$, $360$, $432$, $1080$ and $1080$ singular points.\finl
\end{itemize}
\end{exemple}

\bigskip

The other exceptional groups have been investigated but the computations 
are somewhat too long (note that $n \ge 5$).

\bigskip

\section{The case of ${\boldsymbol{G_{32}}}$}\label{sec:g32}

\medskip

\boitegrise{{\bf Hypothesis.} {\it We assume in this section, and only in this section, 
that $W$ is the primitive complex reflection group $G_{32}$.}}{0.75\textwidth}

\bigskip

In Table~\ref{tab:singular surfaces}, it is said that the surface $\ZC(F_{u_3})$ 
attached to $G_{32}$ has $1~\!440$ singularities of type $D_4$. We give here a detailed account of this example, 
and show that it also produces surfaces of degree $8$ and $16$ with many singularities 
of type $D_4$. The {\sc Magma} codes are contained in the {\sc arXiv} version of 
this section~\cite{1440}. 

\bigskip

We need some more notation. 
If $f \in \CM[x_1,x_2,x_3,x_4]$ is homogeneous, 
we denote by $f[k]$ the homogeneous polynomial 
$f(x_1^k,x_2^k,x_3^k,x_4^k)$. Let $W_1$ be the subgroup of $\Gb\Lb_4(\CM)$ generated by
$$
s_1=\begin{pmatrix}
0 &  1 &  0 &  0 \\
1 &  0 &  0 &  0 \\
0 &  0 &  1 &  0 \\
0 &  0 &  0 & - 1\\
\end{pmatrix},
\quad s_2 =\begin{pmatrix}
1&  0&  0&  0\\
0&  0&  1&  0\\
0&  1&  0&  0\\
0&  0&  0& - 1\\
\end{pmatrix}\quad\text{and}\quad
s_3=\begin{pmatrix}
- 1 &  0 &  0 &  0 \\
0  &  1 &  0 &  0 \\
0  &  0 &  0 &  1 \\
0  &  0 &  1 &  0 \\
\end{pmatrix}.
$$
Let $\z_3$ (resp. $\z_4$) be a primitive third (resp. fourth) root of unity. 
Let $W_2$ be the subgroup of $\Gb\Lb_4(\CM)$ generated by
$$s_1'=\begin{pmatrix}
0 &  1 &  0 &  0\\
1 &  0 &  0 &  0\\
0 &  0 &  1 &  0\\
0 &  0 &  0 &  \z_4\\
\end{pmatrix},
\quad
s_2'=\begin{pmatrix}
1 &  0 &  0 &  0 \\
0 &  0 &  1 &  0 \\
0 &  1 &  0 &  0 \\
0 &  0 &  0 &  \z_4 \\
\end{pmatrix}
\quad\text{and}\quad
s_3'=\begin{pmatrix}
-\z_4&  0 &  0 &  0 \\
   0 &  1 &  0 &  0 \\
   0 &  0 &  0 &  1 \\
   0 &  0 &  1 &  0 \\
\end{pmatrix}.$$
Finally, let $W=W_3$ denote the subgroup of $\Gb\Lb_4(\CM)$ generated by 
$$s_1''=\begin{pmatrix}
1 & 0 & 0    & 0 \\
0 & 1 & 0    & 0 \\
0 & 0 & \z_3 & 0 \\
0 & 0 & 0    & 1 \\
\end{pmatrix},\quad
s_2''=\begin{pmatrix}
\frac{\z_3+2}{3} & \frac{\z_3-1}{3} & \frac{\z_3-1}{3} & 0 \\
\frac{\z_3-1}{3} & \frac{\z_3+2}{3} & \frac{\z_3-1}{3} & 0 \\
\frac{\z_3-1}{3} & \frac{\z_3-1}{3} & \frac{\z_3+2}{3} & 0 \\
0 & 0 & 0 & 1 \\
\end{pmatrix},$$
$$
s_3''=\begin{pmatrix}
1 & 0 & 0 & 0 \\
0 & \z_3 & 0 & 0 \\
0 & 0 & 1 & 0 \\
0 & 0 & 0 & 1 \\
\end{pmatrix}
\quad\text{and}\quad
s_4''=\begin{pmatrix}
\frac{\z_3+2}{3} & \frac{1-\z_3}{3} & 0 & \frac{1-\z_3}{3} \\
\frac{1-\z_3}{3} & \frac{\z_3+2}{3} & 0 & \frac{\z_3-1}{3} \\
0 & 0 & 1 & 0 \\
\frac{1-\z_3}{3} & \frac{\z_3-1}{3} & 0 & \frac{\z_3+2}{3} \\
\end{pmatrix}.$$

\bigskip

\noindent{\bf Commentaries.} 
The following facts are checked using {\sc Magma}, as explained in~\cite{1440}. 
Let $\Zrm(W_i)$ denote the center of $W_i$. In all cases, 
it is isomorphic to a group of roots of unity acting by scalar 
multiplication. Then:
\begin{itemize}
\itemth{a} The group $W_1$ has order $48$ and is isomorphic to 
the non-trivial double cover $\Gb\Lb_2(\FM_{\! 3})$ of the symmetric group 
$\SG_4 \simeq W_1/\Zrm(W_1)$.

\itemth{b} The group $W_2$ has order $768$, 
contains a normal abelian subgroup $H$ of order $32$ 
and $W_2/H \simeq \SG_4$. The group $W_2/\Zrm(W_2)$ has order 
$192$, but is not isomorphic 
to a Coxeter group of type $D_4$.

\itemth{c} The group $W_3$ is the complex reflection group denoted 
by $G_{32}$ in the Shephard--Todd classification~\cite{shephard todd} (it has order $155~\!920$). 
Recall that the group $W_3/\Zrm(W_3)$ is a simple group of order $25~\!920$ 
and is isomorphic to the derived subgroup of the Weyl group of type $E_6$ (i.e. 
to the derived subgroup of the special orthogonal group 
$\Sb\Ob_5(\FM_{\! 3})$). It contains the group $W_1$ as a subgroup, as well as 
a subgroup of diagonal matrices isomorphic to $(\mub_3)^4$, where $\mub_d$ 
is the group of $d$-th roots of unity.

Note that we have used the version of $G_{32}$ implemented by Michel in the 
{\tt Chevie} package of {\tt GAP3}~\cite{jean}.\finl
\end{itemize}

\bigskip

If $\l=(\l_1 \ge \l_2 \ge \l_3 \ge \l_4)$ is a partition of $8$ 
of length at most $4$, we denote by $\O_\l^-$ (resp. $\O_\l^+$) be the orbit of the monomial 
$x_1^{\l_1}x_2^{\l_2}x_3^{\l_3}x_4^{\l_4}$ under the action of $W_1$ (resp. the symmetric group 
$\SG_4$) and we set
$$m_\l^\e=\sum_{m \in \O_\l^\e} m$$
for $\e \in \{+,-\}$. 
Then $m_\l^+$ is the symmetric function traditionnally 
denoted by $m_\l$. If all the $\l_i$'s are even, then 
$m_\l^-=m_\l^+$ but note for instance that
\eqna
m_{611}^+ \neq m_{611}^- &=&  
x_1^6 x_2 x_3 + x_1^6 x_2 x_4 - x_1^6 x_3 x_4 + x_1 x_2^6 x_3 - x_1 x_2^6 x_4 + x_2^6 x_3 x_4 \\ 
&& + x_1 x_2 x_3^6 + x_1 x_3^6 x_4 - x_2 x_3^6 x_4 - x_1 x_2 x_4^6 - x_1 x_3 x_4^6 - x_2 x_3 x_4^6. 
\endeqna
Now, let
\eqna
g\!\!\!\!&=&\!\!\!\!m_8^- - 6 m_{62}^- - 60 m_{611}^- + 2~\!240 m_{521}^- -14 m_{44}^- 
+ 10~\!180 m_{431}^- + 40~\!412 m_{422}^- \\ && - 23~\!440 m_{4211}^- 
+ 111~\!980 m_{332}^- + 154~\!704 m_{2222}^-.
\endeqna
By construction, $m_\l^-$ is invariant under the action of $W_1$ and so 
$g$ is invariant under the action of $W_1 \simeq \tilde{\SG}_4$. One can check with {\sc Magma} 
the following facts~\cite[Proposition~1]{1440}:

\bigskip

\begin{prop}\label{prop 1}
If $1 \le k \le 3$, then the polynomial $g[k]$ 
is invariant under the action of $W_k$.
\end{prop}

\bigskip

One can also check that $g[3]$ is the polynomial denoted by $F_{u_3}$ (suitably normalized) 
in Table~\ref{tab:singular surfaces} (in the $G_{32}$ example). 

\bigskip

\begin{theo}\label{theo 2}
The homogeneous polynomial $g$ satisfies the following statements:
\begin{itemize}
\itemth{a} $\ZC(g)$ is an irreducible surface of degree $8$ in $\Pb^3(\CM)$ 
with exactly $44$ singular points which are all quotient singularities of 
type $D_4$.

\smallskip

\itemth{b} If $k \ge 1$, then $\ZC(g[k])$ is an irreducible surface 
of degree $8k$, whose singular locus has dimension $0$ and contains at least $44k^3$ quotient 
singularities of type $D_4$.

\smallskip

\itemth{c} $\ZC(g[2])$ is an irreducible surface of degree $16$ with exactly $472$ singular points: 
$24$ quotient singularities of type $A_1$, $96$ quotient singularities of type $A_2$ and 
$352$ quotient singularities of type $D_4$.

\smallskip

\itemth{d} $\ZC(g[3])$ is an irreducible surface of degree $24$ in $\Pb^3(\CM)$ 
with exactly $1~\!440$ singular points which are all quotient singularities of type $D_4$. 
\end{itemize}
\end{theo}

\bigskip

\begin{rema}
Note that $g$ has coefficients in $\QM$ but the singular points 
of $\ZC(g)$, $\ZC(g[2])$ and $\ZC(g[3])$ have coordinates 
in various field extensions of $\QM$, and most of the singular points 
are not real (at least in this model).\finl
\end{rema}

\bigskip

We now turn to the study of the singularities 
of the varieties $\ZC(g[i])$ for $i \in \{1,2,3\}$. Note the 
following fact, checked using {\sc Magma}~\cite[Lemma~3]{1440}, that will be used further:

\bigskip

\begin{lem}\label{lem 3}
If $1 \le i < j \le 4$, then the closed 
subscheme of $\Pb^3(\CM)$ defined by the homogeneous ideal 
$\langle g,\frac{\partial g}{\partial x_i},\frac{\partial g}{\partial x_j} \rangle$ 
has dimension $0$.
\end{lem}

\bigskip

\subsection{Degree ${\boldsymbol{8}}$}\label{sec:8}
The {\sc Magma} computations leading to the 
proof of the statement~(a) of Theorem~\ref{theo 2} are detailed in~\cite[\S{1}]{1440}. 
Along these computations, the following facts are obtained 
(here, $\UC$ denotes the open subset of $\Pb^3(\CM)$ defined by $x_1x_2x_3x_4 \neq 0$):

\bigskip

\begin{prop}\label{prop 4}
We have:
\begin{itemize}
\itemth{a} $\dim \ZC_\sing(g) = 0$, so $\ZC(g)$ is irreducible. 

\itemth{b} $\ZC_\sing(g)$ is contained in $\UC$.

\itemth{c} The group $W_1$ has $3$ orbits in $\ZC_\sing(g)$, of respective length $8$, $12$ and $24$. 
\end{itemize}
\end{prop}

\bigskip

\begin{center}
\begin{figure}
\begin{center}
\includegraphics[scale=0.2]{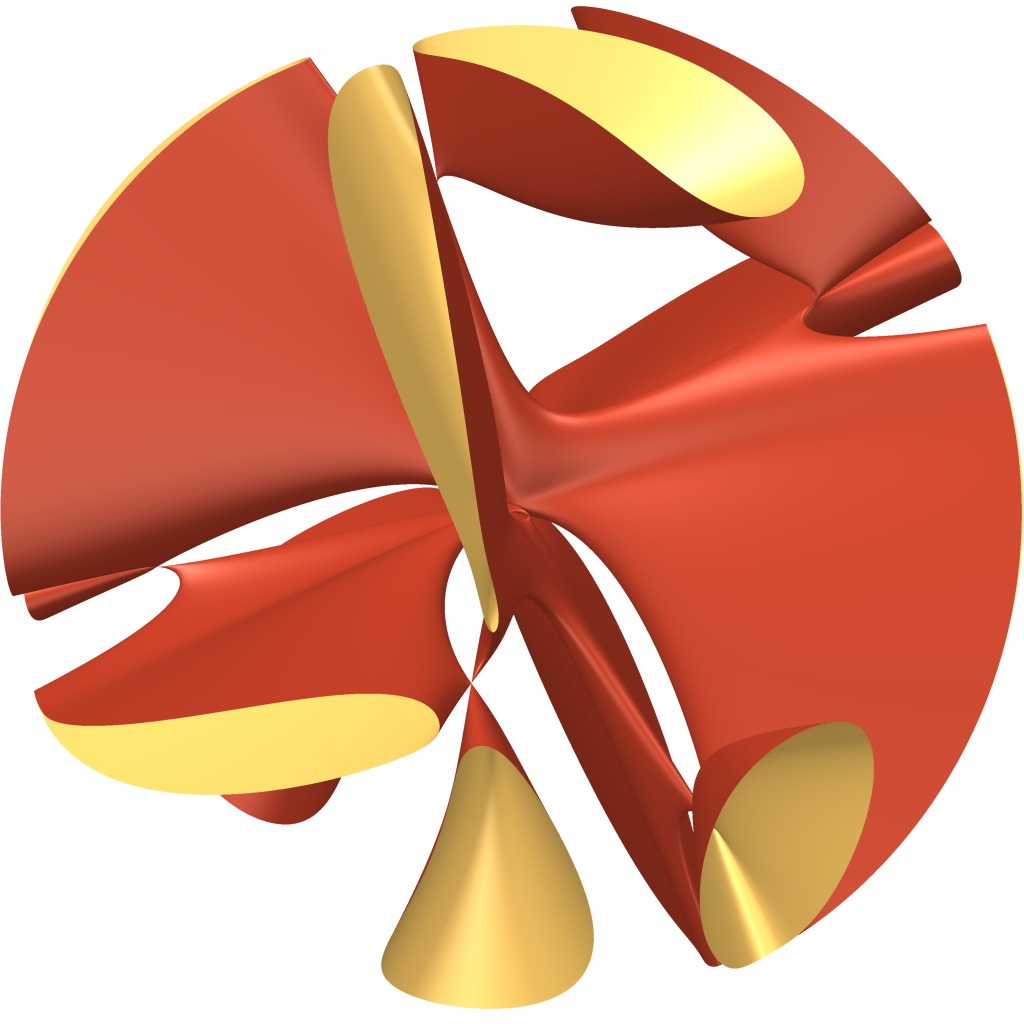}
\includegraphics[scale=0.2]{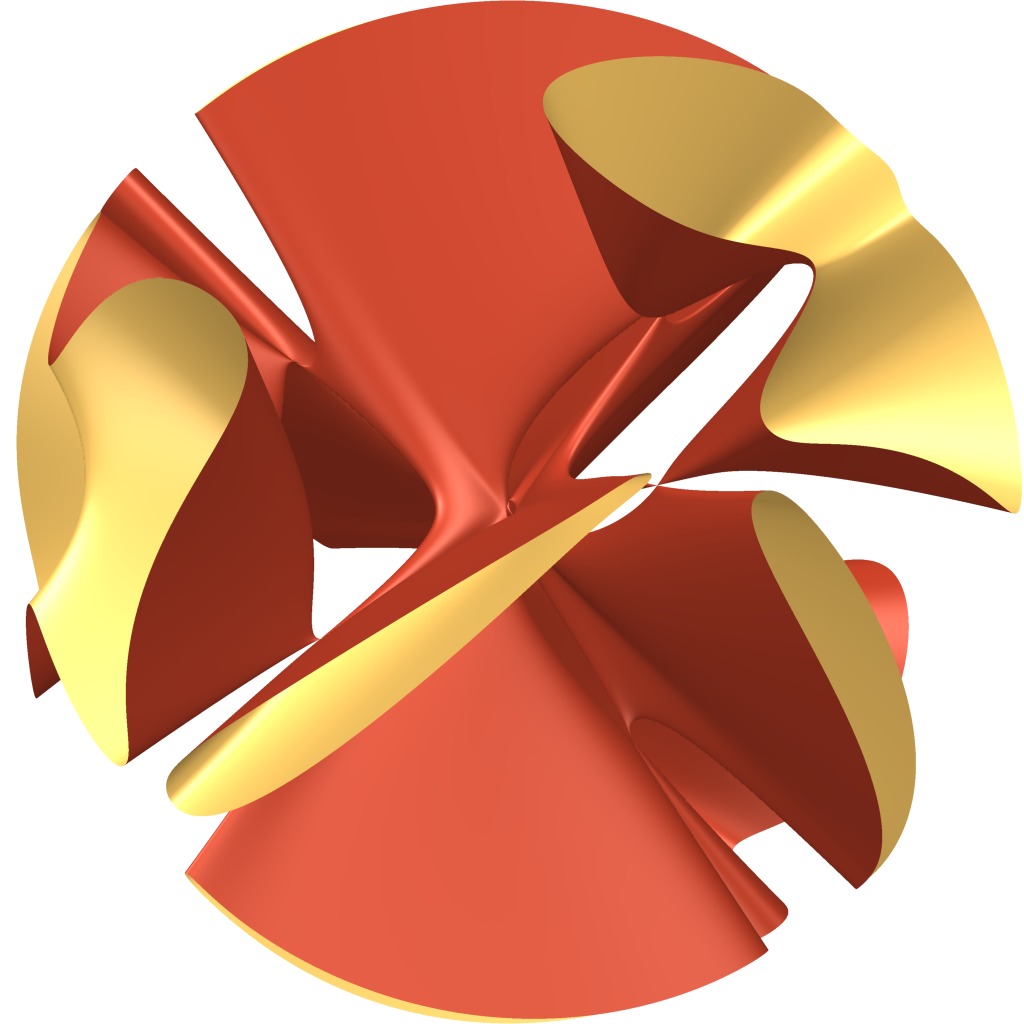}
\caption{Part of the real locus of $\ZC(g)$}\label{fig:8}
\end{center}
\end{figure}
\end{center}

Note that the points in the $W_1$-orbit of cardinality $8$ 
are the only real singular points of $\ZC(g)$. 
Figure~\ref{fig:8} shows part of the real locus of $\ZC(g)$. 

\bigskip

\subsection{Degree ${\boldsymbol{8k}}$}\label{sec:8k}
Let $\UC$ denote the open subset of $\Pb^3(\CM)$ defined by $x_1x_2x_3x_4 \neq 0$ 
and let $\s_{\! k} : \Pb^3(\CM) \to \Pb^3(\CM)$, $[x_1;x_2;x_3;x_4] \mapsto [x_1^k;x_2^k;x_3^k;x_4^k]$. 
The restriction of $\s_{\! k}$ to a morphism $\UC \to \UC$ is an \'etale Galois covering, 
with group $(\mub_k)^4/\D\mub_k$ (here, $\D : \mub_k \injto (\mub_k)^4$ 
is the diagonal embedding). We have $\ZC(g[k])=\s_{\! k}^{-1} (\ZC(g))$. 

Let us first prove that $\ZC(g[k])$ is irreducible. We may assume that 
$k \ge 2$, as the result has been proved for $k=1$ in the previous section. 
Recall that
$$\frac{\partial g[k]}{\partial x_i}=kx_i^{k-1}(\frac{\partial g}{\partial x_i} \circ \s_{\! k}),$$
so the singular locus of $\ZC(g[k])$ is contained in 
$$\{p_1,p_2,p_3,p_4\} \cup \Bigl(\bigcup_{i \neq j} \s_{\! k}^{-1} (\ZC_{i,j})\Bigr),$$
where $p_i=[\d_{i1};\d_{i2};\d_{i3};\d_{i4}]$ (and $\d_{ij}$ is the Kronecker symbol) and 
$\ZC_{i,j}$ is the subscheme of $\Pb^3(\CM)$ defined by the ideal 
$\langle g,\frac{\partial g}{\partial x_i},\frac{\partial g}{\partial x_j} \rangle$ 
(and which has dimension $0$ by Lemma~\ref{lem 3}). 
Since $\s_{\! k}$ is finite, this implies that $\ZC_\sing(g[k])$ has dimension $0$, 
so $\ZC(g[k])$ is irreducible.

Now, $\s_{\! k} : \UC \to \UC$ is \'etale and the singular locus of $\ZC(g)$ 
is contained in $\UC$ (see Proposition~4(b)). Therefore, the $44$ singularities of $\ZC(g)$ lift to  
$44 k^3$ singularities in $\ZC(g[k]) \cap \UC$ of the same type, i.e. quotient 
singularities of type $D_4$. This proves the statement~(b) of Theorem~\ref{theo 2}.

Note that, for $k=2$, $3$ and $4$ (and maybe for bigger $k$) we will prove 
in the next sections that $\ZC(g[k])$ contains singular points outside 
of $\UC$. 

\bigskip

\subsection{Degree ${\boldsymbol{16}}$}\label{sec:16}
\begin{center}
\begin{figure}
\begin{center}
\includegraphics[scale=0.2]{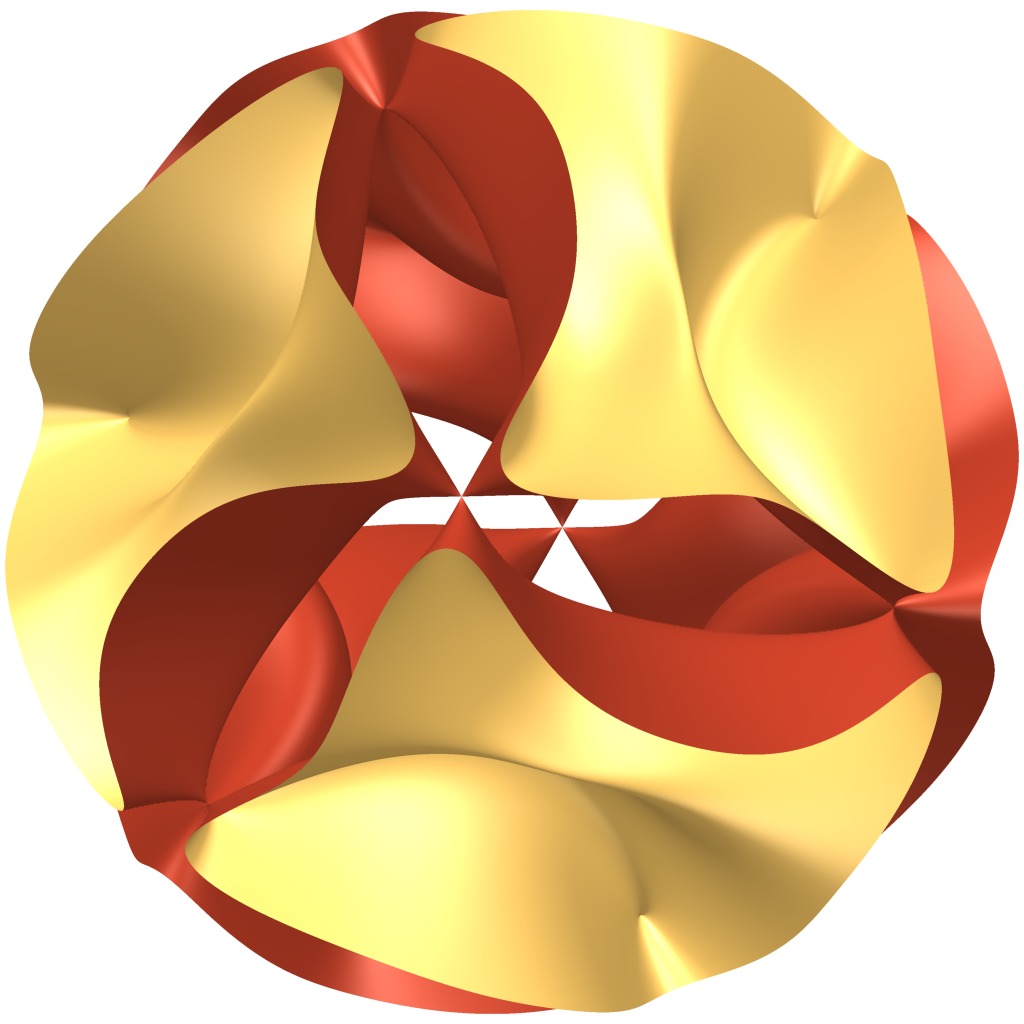}
\includegraphics[scale=0.2]{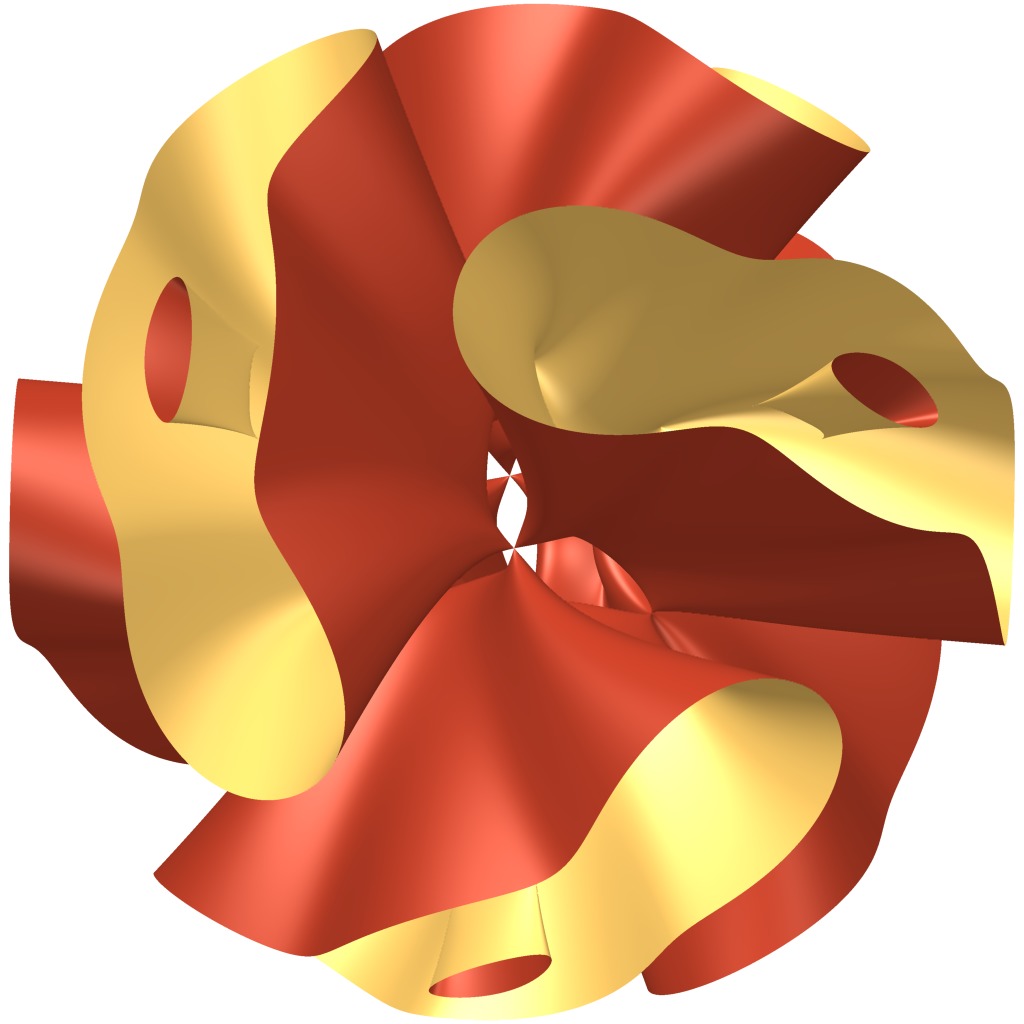}
\caption{Part of the real locus of $\ZC(g[2])$}\label{fig:16}
\end{center}
\end{figure}
\end{center}

Using the morphism $\s_{\! 2}$ defined in the previous section, we get that $\ZC(g[2]) \cap \UC$ 
has exactly $352$ singular points, which are all quotient singularities 
of type $D_4$. The other singularities are determined thanks to {\sc Magma} 
computations that are detailed in~\cite[\S{3}]{1440}, 
and which confirm the statement~(c) of Theorem~\ref{theo 2}. Note that we also 
need the software {\sc Singular}~\cite{singular} for computing some 
Milnor numbers and identifying the singularity $A_2$. Note also that $W_2$ acts transitively on the 
$24$ quotient singularities of type $A_1$ and also acts transitively on the 
$96$ quotient singularities of type $A_2$. 
Figure~\ref{fig:16} shows part of the real locus of $\ZC(g)$.

\subsection{Degree ${\boldsymbol{24}}$}\label{sec:24}
\begin{figure}
\includegraphics[scale=0.2]{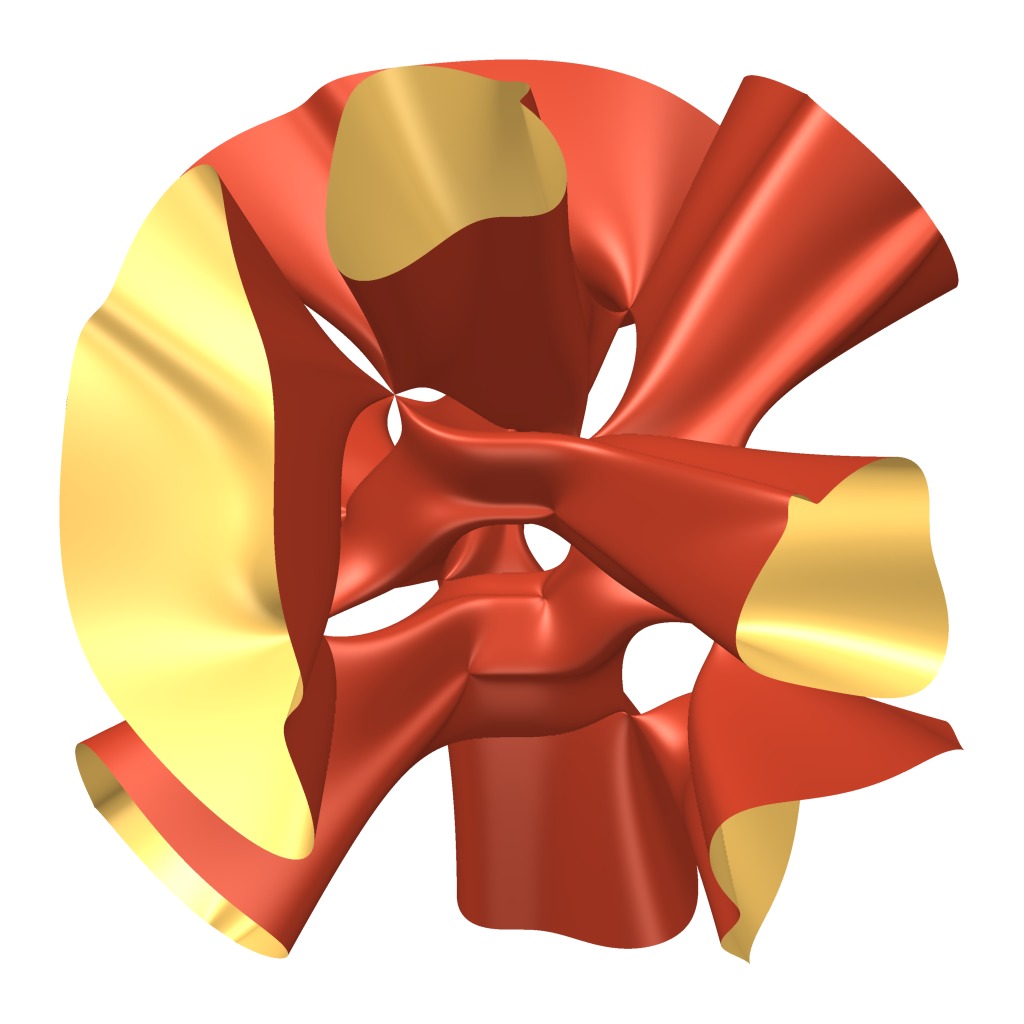}
\includegraphics[scale=0.2]{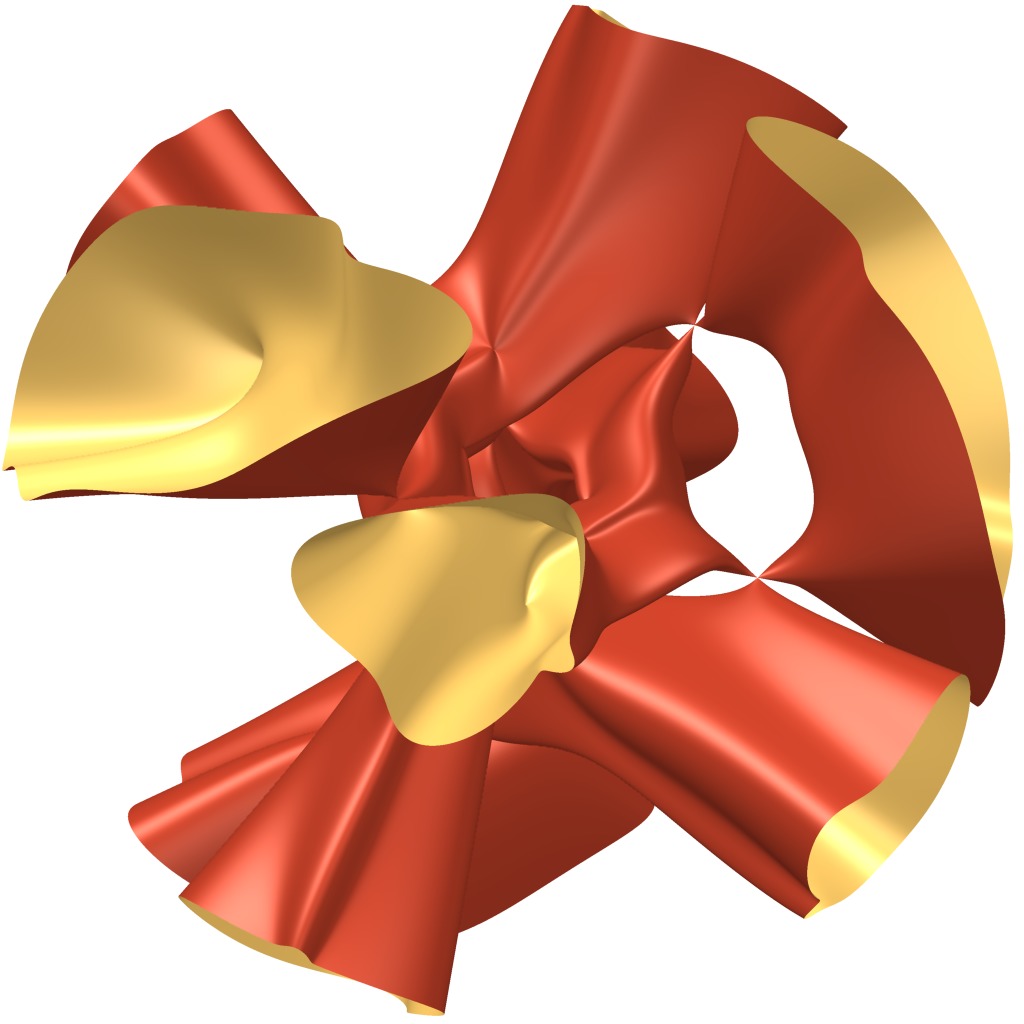}
\caption{Part of the real locus of $\ZC(g[3])$}\label{fig:24}
\end{figure}

Using the morphism $\s_{\! 3}$ defined in Section~\ref{sec:8k}, we get that $\ZC(g[3]) \cap \UC$ 
has exactly $44 \times 3^3=1~\!188$ singular points, which are all quotient singularities 
of type $D_4$. The other singularities are determined thanks to {\sc Magma} 
computations that are detailed in~\cite{bonnafe calculs} or~\cite[\S{4}]{1440}, 
and which confirm the statement~(d) of Theorem~\ref{theo 2}. Note also that, in the given model, 
the surface $\ZC(g[3])$ has only $32$ real singular points: Figure~\ref{fig:24} gives partial views 
of its real locus.

\bigskip

\subsection{Complements}
From Section~\ref{sec:8k}, we deduce that $\ZC_\sing(g[4])$ has $2~\!816$ quotient 
singularities of type $D_4$ lying in the open subset $\UC$ and it can be checked that 
it has $432$ other singular points not lying in $\UC$, for which we did not determine
the type.

\bigskip

\noindent{\bf Acknowledgements.} This paper is based 
upon work supported by the National Science Foundation under 
Grant No. DMS-1440140 while the author was in residence at the Mathematical Sciences Research 
Institute in Berkeley, California, during the Spring 2018 semester. 
The hidden computations which led to the discovery 
of the polynomial $g$ were done using the High Performance Computing facilities 
of the MSRI.

I wish to thank warmly Alessandra Sarti, Oliver Labs and Duco van Straten 
for useful comments and references and Gunter Malle for a careful reading 
of a first version of this paper. 
Figures were realized using the software {\sc SURFER}~\cite{surfer}.


\bigskip


\begin{thebibliography}{AAAAa}
\bibitem[AGV]{agv} {\sc V. I. Arnold, S. M. Gusein-Zade \& A. N. Varchenko}, 
{\it Singularities of
differentiable maps}, Vol. I, The classification of critical points, caustics and
wave fronts, volume {\bf 82} of Monographs in Mathematics. Birkh\"auser Boston
Inc., Boston, MA, 1985.

\bibitem[Bar]{barth} {\sc W. Barth}, 
{\it Two projective surfaces with many nodes, admitting the symmetries of the icosahedron}, 
J. Algebraic Geom. {\bf 5} (1996), 173--186.

\bibitem[Bon1]{1440} {\sc C. Bonnaf\'e}, 
{\it A surface of degree $24$ with $1440$ singularities of type $D_4$}, 
preprint (2018), {\tt arXiv:1804.08388}.

\bibitem[Bon2]{bonnafe calculs} {\sc C. Bonnaf\'e}, 
{\it Magma codes for ``Some singular curves and surfaces arising from invariants 
of complex reflection groups''}, preprint (2018), available at 
{\tt hal.archives-ouvertes.fr/hal-01897587}; also available at 
{\tt http://imag.umontpellier.fr/\~{}bonnafe/}.

\bibitem[Bro]{broue} {\sc M. Brou\'e}, 
{\it Introduction to complex reflection groups and their braid groups}, 
Lecture Notes in Mathematics {\bf 1988}, Springer-Verlag, Berlin, 2010, xii+138 pp.

\bibitem[Bru]{bruce} {\sc J.W. Bruce}, 
{\it An upper bound for the number of singularities on a projective hypersurface}, 
Bull. Lond. Math. Soc. {\bf 13} (1981), 47-51.

\bibitem[Chm]{chmutov} {\sc S. V. Chmutov}, 
{\it Examples of projective surfaces with many singularities}, 
J. Algebraic Geom. {\bf 1} (1992), 191-196.

\bibitem[C-ALi]{39} {\sc J.-I. Cogolludo-Agust\'{\i}n \& A. Libgober}, 
{\it Mordell-Weil groups of elliptic threefolds and the Alexander module of plane curves}, 
J. Reine Angew. Math. {\bf 697} (2014), 15-55.

\bibitem[DGPS]{singular}
{\sc W. Decker, G.-M. Greuel, G. Pfister \& H. Sch{\"o}nemann},  
{\it Singular {4-1-1} --- {A} computer algebra system for polynomial computations}, 
{\tt http://www.singular.uni-kl.de} (2018).

\bibitem[End]{endrass} {\sc S. Endra}\ss, 
{\it A projective surface of degree eight with 168 nodes}, J. Algebraic Geom. {\bf 6} (1997), 
325-334.

\bibitem[EnPeSt]{triple} {\sc S. Endra}\ss, {\sc U. Persson \& J. Stevens}, 
{\it Surfaces with triple points}, 
J. Algebraic Geom. {\bf 12} (2003), 367-404.

\bibitem[Esc1]{escudero} {\sc J. G. Escudero}, 
{\it Hypersurfaces with many $A_j$ singularities: explicit constructions}, 
J. Comput. Appl. Math. {\bf 259} (2014), 87-94.

\bibitem[Esc2]{escudero 2} {\sc J. G. Escudero}  
{\it Arrangements of real lines and surfaces with $A$ and $D$ singularities}, 
Experimental Mathematics {\bf 23} (2014), 482-491.

\bibitem[Iv]{ivinskis} {\sc K. Ivinskis}, 
{\it Normale Fl\"achen und die Miyaoka-Kobayashi Ungleichung}, 
Diplomarbeit (Bonn), 1985.

\bibitem[Lab]{labs} {\sc O. Labs}, 
{\it A septic with 99 real nodes}, 
Rend. Sem. Mat. Univ. Pad. {\bf 116} (2006), 299-313.
 
\bibitem[Magma]{magma} {\sc W. Bosma, J. Cannon \& C. Playoust}, 
{\it The Magma algebra system. I. The user language}, J. Symbolic Comput. {\bf 24} (1997), 235-265. 

\bibitem[MaMi]{marin-michel} {\sc I. Marin \& J. Michel}, 
{\it Automorphisms of complex reflection groups}, Represent. Theory {\bf 14} (2010), 747-788.

\bibitem[Mic]{jean} {\sc J. Michel}, 
{\it The development version of the CHEVIE package of GAP3}, 
J. of Algebra {\bf 435} (2015), 308--336.

\bibitem[Miy]{miyaoka} {\sc Y. Miyaoka}, 
{\it The maximal number of quotient singularities on surfaces with given numerical invariants}, 
Math. Ann. {\bf 268} (1984), 159-171.

\bibitem[Sak]{sakai} {\sc F. Sakai}, {\it Singularities of plane curves}, 
Geometry of complex projective varieties (Cetraro, 1990), 257-273,
Sem. Conf., {\bf 9}, Mediterranean, Rende, 1993.

\bibitem[Sar1]{sarti 0} {\sc A. Sarti}, 
{\it Pencils of Symmetric Surfaces in $\Pb^3$}, 
J. Algebra {\bf 246} (2001), 429-452.

\bibitem[Sar2]{sarti 1} {\sc A. Sarti},
{\it Symmetric surfaces with many singularities}, 
Comm.  Algebra {\bf 32} (2004), 3745-3770.

\bibitem[Sar3]{sarti} {\sc A. Sarti}, 
{\it Symmetrische Fl\"{a}chen mit gew\"{o}hnlichen Doppelpunkten}, 
Math. Semesterber. {\bf 55} (2008), 1-5. 

\bibitem[ShTo]{shephard todd} {\sc G. C. Shephard \& J. A. Todd}, 
{\it Finite unitary reflection groups}, 
Canad. J. Math. {\bf 6} (1954), 274-304.

\bibitem[Sta]{stagnaro} {\sc E. Stagnaro}, 
{\it A degree nine surface with 39 triple points}, 
Ann. Univ. Ferrara Sez. VII {\bf 50} (2004), 111-121. 

\bibitem[Sur]{surfer} {\tt www.imaginary.org/program/surfer}.

\bibitem[Var]{varchenko} {\sc A.N. Varchenko}, 
{\it On the semicontinuity of the spectrum and an upper bound for the number
of singular points of a projective hypersurface}, 
J. Soviet Math. {\bf 270} (1983), 735-739.

\bibitem[Wah]{wahl} {\sc J. Wahl}, 
{\it Miyaoka-Yau inequality for normal surfaces and local analogues}, in 
Classification of algebraic varieties (L'Aquila, 1992), 381-402,
Contemp. Math. {\bf 162}, Amer. Math. Soc., Providence, RI, 1994. 

\bibitem[Thi]{thiel} {\sc U. Thiel}, 
{\it CHAMP: A Cherednik Algebra Magma Package}, 
LMS J. Comput. Math. {\bf 18} (2015), 266-307. 
\end{thebibliography}
\end{document}